\documentclass[reqno,centertags,a4paper,11pt]{amsart}
\usepackage[a4paper,margin=2.54cm]{geometry}

\usepackage{dsfont}
\usepackage{amscd}
\usepackage{amsmath,tikz}
\usepackage{amssymb}
\usepackage{mathtools}
\usepackage{amsfonts}
\usepackage{enumerate}
\usepackage{latexsym}
\usepackage{cases}
\usepackage{amsthm}
\usepackage{graphicx}
\usepackage{verbatim}
\usepackage{mathrsfs}
\usepackage{color}
\usepackage{wasysym}
\usepackage{makecell}
\usepackage{enumitem}
\usepackage{tikz}
\usepackage{subcaption}
\usepackage[numbers]{natbib}
\usepackage{hyperref}

\newcommand{\on}{\operatorname}

\newtheorem{theorem}{Theorem}[section]
\newtheorem{lemma}{Lemma}[section]

\newtheorem{definition}[lemma]{Definition}

\newtheorem{remark}[lemma]{Remark}
\theoremstyle{definition}
\numberwithin{equation}{section}

\newcommand{\dd}{\mathrm{d}}

\newcommand{\R}{\mathbb{R}}
\newcommand{\N}{\mathbb{N}}

\newcommand{\supp}{\mathrm{supp}}

\allowdisplaybreaks[2]

\title[restriction estimates for degenerate higher codimensional quadratic surfaces]{Sharp restriction estimates for some degenerate higher codimensional quadratic surfaces}
\author{Zhenbin Cao\textsuperscript{1}}
\thanks{\textsuperscript{1}Institute of Mathematics, Henan Academy of Sciences, Zhengzhou 450046, China}

\author{Changxing Miao\textsuperscript{2}}
\thanks{\textsuperscript{2}School of Mathematics and Physics, University of Science and Technology Beijing, Beijing 100083, China}

\author{Yixuan Pang\textsuperscript{3,*}}
\thanks{\textsuperscript{3}Department of Mathematics, University of Pennsylvania, Philadelphia, PA 19104}
\thanks{\textsuperscript{*}Correspondence to be sent to: pyixuan@sas.upenn.edu}


\begin{document}

\begin{abstract}
The Fourier restriction conjecture is a fundamental problem in harmonic analysis. In this paper, we investigate restriction estimates for degenerate higher codimensional quadratic surfaces and obtain sharp results for some types of degenerate cases. A major obstacle in establishing sharp restriction estimates is the failure of rescaling invariance, which is crucial for induction on scale to be effective. Motivated by the work of Guo and Oh (2022), we introduce a method, building on an iterative variant of the broad--narrow analysis, that does not heavily rely on induction on scale. To obtain suitable transversality conditions for this analysis and to derive desirable bounds for the broad part, we define a generalized notion of Jacobian, and establish its structural properties. These properties are proved using tools and techniques from both algebra and graph theory. 
\end{abstract}

\maketitle

\section{Introduction}
Let $d,n \geq 1$, and $\mathbf{Q}(\xi)=(Q_1(\xi),...,Q_n(\xi))$ be an $n$-tuple of real quadratic forms on $\R^d$. The graph of such a tuple, $S_{\mathbf{Q}}=\{ (\xi,\mathbf{Q}(\xi))\in [0,1]^d \times \mathbb{R}^n \}$, is a $d$-dimensional submanifold of $\mathbb{R}^{d+n}$. Define the extension operator associated to $S_{\mathbf{Q}}$ by
$$E^{\mathbf{Q}}f(x):=\int_{[0,1]^d} e^{ 2\pi ix\cdot(\xi, \mathbf{Q}(\xi))}f(\xi)\dd\xi,\quad \quad x\in \mathbb{R}^{d+n}.$$
We mainly focus on the following Fourier restriction problem: Find optimal ranges of $p$ and $q$ such that
\begin{equation}\label{main aim}
    \|E^{\mathbf{Q}}f\|_{L^q(\mathbb{R}^{d+n})}\leq C_{p,q,\mathbf{Q},d,n}\|f\|_{L^p(\mathbb{R}^{d})} , 
\end{equation}
for every measurable function $f$.

Originating from a deep observation of Stein \cite{stein1979some} in 1967, the Fourier restriction conjecture has been widely studied for hypersurfaces ($n=1$), especially the paraboloid. Although the $d=1$ case was solved by Fefferman \cite{Fef70} and Zygmund \cite{Zyg74} half a century ago, the $d\geq 2$ case is still far from being fully understood. In the special case when $p=2$, the sharp estimate up to the endpoint for the sphere was first proved by Tomas \cite{Tomas75}, and the endpoint result was later established by Stein \cite{stein1986osc} through complex interpolation. The Stein--Tomas framework is very influential and turns out to work for any hypersurface with nonvanishing Gaussian curvature, including the paraboloid. However, things become much trickier for general $p$, and people have to search for other methods. After early major progress made by Bourgain \cite{Bou91}, several influential methods have been developed, such as  the bilinear method \cite{Tao98,Tao03}, the multilinear method/broad--narrow analysis \cite{BG11}, and the polynomial method \cite{Gut16,Gut18}. These methods have been combined with more delicate techniques from incidence geometry and real algebraic geometry \cite{HR19,HZ23,Wan22,GWZ24}. Recently, Wang--Wu \cite{WW24} set the current world records for $d = 2$ and large $d$ by developing a novel framework which combines refined decoupling with deep tools from geometric measure theory and additive combinatorics.

In contrast to the paraboloid, people know much less about restriction estimates for higher codimensional quadratic surfaces ($n>1$). We now briefly summarize the history.

For a long time, the state-of-the-art methods are hardly beyond the Stein--Tomas framework and interpolation of analytic families of operators, so the $p=2$ case is the main focus of many papers. When $n=2$, Christ \cite{Chr82} and Mockenhaupt \cite{Moc96} introduced a geometric notion of nondegeneracy which is equivalent to optimal $L^2\rightarrow L^q$ estimates. However, when $n\geq 3$, things become much more complicated and existing research is fragmented, see \cite{banner2002restriction} for some partial progress. 

Another line of attack is based on the Fefferman--Zygmund framework and mapping properties of multilinear fractional integral operators. Such a method breaks through the $L^2$ barrier, and can sometimes yield full sharp ranges of $L^p\rightarrow L^q$ estimates (even including the critical line). For example, when $d=n=2$, Christ \cite{Chr85} fully solved the restriction problem for non-degenerate quadratic surfaces, and when $n=\frac{d(d+1)}{2}$, Bak--Lee \cite{BL04} and Oberlin \cite{oberlin2005restriction} fully solved the restriction problem for extremal quadratic surfaces. However, such a method seems to require certain unnatural relations between $d$ and $n$, which prevents us from obtaining sharp results in general cases, see \cite{oberlin2004convolution} for a discussion on the non-degenerate surface $(\xi_1,\xi_2,\xi_3,\xi_1^2+\xi_2^2, \xi_2^2+\xi_3^2)$ using similar techniques.

In recent years, the developments when $n=1$ mentioned before have greatly facilitated the research when $n>1$. For example, Bak--Lee--Lee \cite{BLL17} applied the bilinear method to obtain some restriction estimates for higher codimensional surfaces with nonvanishing rotational curvature, such as the complex paraboloid $(\xi_1,\xi_2,\xi_3,\xi_4,\xi_1^2-\xi_2^2+\xi_3^2-\xi_4^2, \xi_1\xi_2 + \xi_3\xi_4)$, and Lee--Lee \cite{LL20} applied the broad--narrow analysis to study holomorphic complex hypersurfaces. However, the assumptions in these two works are too strong due to technical reasons. Noteworthy progress was made by Guo--Oh \cite{GO22}, who systematically investigated restriction estimates for all quadratic surfaces of codimension 2 in $\mathbb{R}^5$. Later on, Guo--Oh--Zhang--Zorin-Kranich \cite{guo2023decoupling} established sharp decoupling inequalities for all higher codimensional quadratic surfaces by combining a transversality condition originated from \cite{bourgain2016mean}. As an application, they obtained a unified restriction estimate for quadratic surfaces of any dimensions and codimensions through the broad--narrow analysis. However, the full power of the broad--narrow analysis has not been exploited until Gan--Guth--Oh \cite{GGO23} further developed the framework in \cite{guo2023decoupling} by adopting a more delicate notion of transversality, which enables them to capture lower dimensional contribution. 

In this paper, we consider an iterative version of the broad--narrow analysis from \cite{GGO23,guo2023decoupling}. By excavating structural property sufficiently for each individual quadratic form, we can build sharp restriction estimates for several classes of degenerate higher codimensional quadratic surfaces. We take some representative examples as below. Before that, we first introduce a generalized notion of the Jacobian: For a tuple $\mathbf{Q}(\xi)=(Q_1(\xi),...,Q_n(\xi))$ of quadratic forms on $\mathbb{R}^d$ and $d\geq n$, we define 
\begin{equation*}
    J_{\mathbf{Q}}(\xi;i_1,...,i_{d-n}):= \begin{vmatrix}
        \nabla \mathbf{Q}   \\
        e_{i_1}  \\
        \vdots\\
        e_{i_{d-n}}
    \end{vmatrix},
\end{equation*}
where $\nabla \mathbf{Q}=(\partial_j Q_i)_{n\times d}$ is a matrix of size $n\times d$, and $e_j$ denotes the $j$th unit coordinate vector in $\mathbb{R}^d$. We adopt the convention that if $d=n$, then $J_\mathbf{Q}(\xi;i_1,...,i_{d-n}) := |\nabla \textbf{Q}(\xi)|$.

\begin{theorem}\label{main thm}
Let $\xi = (\xi_j)_{j=1}^d \in \mathbb{R}^{d}$, and $\tilde{\mathbf{Q}}(\xi) = (\xi_{\lambda_{2m-1}} \xi_{\lambda_{2m}})_{m=1}^n$ be an $n$-tuple of quadratic monomials depending on $\tilde{d}$ variables {\rm(}$\tilde{d} \leq d${\rm)}. Let $P(\xi)$ be a quadratic polynomial, and $\mathbf{Q}$ be the $n$-tuple of quadratic forms obtained by adding $P(\xi)$ to the last component of $\tilde{\mathbf{Q}}$, that is, $\mathbf{Q}(\xi)=\tilde{\mathbf{Q}}(\xi)+(0,...,0,P(\xi))$. Let $w_j$ {\rm(}$1\leq j\leq d${\rm)} be the number of occurrences of $\xi_j$ in $\{ \xi_{\lambda_j}\}_{j=1}^{2n}$, and assume $\max_j w_j\neq1$. Suppose $\tilde{d}\geq n$, $J_{\tilde{\mathbf{Q}}}(\xi;i_1,...,i_{\tilde{d}-n}) \not\equiv 0$ for some $i_1,...,i_{\tilde{d}-n}$, and one of the following two cases holds.

\noindent $\rm{(1)}$ $P(\xi)\equiv0$, and $n\geq 1$.

\noindent $\rm{(2)}$ $P(\xi)\not\equiv0$ is independent of $\{ \xi_{\lambda_j}\}_{j=1}^{2n}$, and $ n\geq2$.

\noindent Then the restriction estimate ${\rm(\ref{main aim})}$ for $\mathbf{Q}$ holds whenever
   \begin{equation}\label{main thm eq2}
   	q > \max_j w_j +2,\quad \frac{1}{p}+\frac{\max_j w_j+1}{q}<1,
   \end{equation}
   which is sharp up to the endpoint for Case $\rm{(1)}$ and, if $\lambda_{2n-1},\lambda_{2n}\notin \{ k: \max_j w_j= w_k \}$, for Case $\rm{(2)}$.
\end{theorem}

\begin{remark}
We briefly comment on the conditions of Theorem ${\rm\ref{main thm}}$. Some of them are crucial. If $\max_j w_j=1$ instead, then $\mathbf{Q}$ must be a tensor product of $n$ hyperbolic paraboloids, whose sharp restriction estimate is far out of reach even in $\R^3$. The constraints ``$\tilde{d}\geq n$" and ``$J_{\tilde{\mathbf{Q}}} \not\equiv 0$'' also play an important role in the proof. In this sense, Case ${\rm(1)}$ is the best one could hope for by our method. Nevertheless, additional requirements in Case ${\rm(2)}$ are not quite essential, but mainly for the convenience and conciseness of the statement of our main theorem. The independence assumption on $P(\xi)$ is to ensure that the polynomial part does not affect the previous monomial parts too much. If there is no restriction ``$\lambda_{2n-1},\lambda_{2n}\notin \{ k: \max_j w_j= w_k \}$,'' then the optimal range for the restriction estimate ${\rm (\ref{main aim})}$ may {\rm(}but not necessarily{\rm)} be different from ${\rm(\ref{main thm eq2})}$. We will introduce more examples  beyond Theorem ${\rm\ref{main thm}}$ in Section~${\rm \ref{section5}}$ {\rm(}see Theorem~{\rm\ref{s5t1})}. The proof of Case ${\rm(2)}$ relies heavily on Case ${\rm(1)}$. We also point out that ``$J_{\tilde{\mathbf{Q}}} \not\equiv 0$'' readily implies ``$J_{\mathbf{Q}} \not\equiv 0$''.
\end{remark}

We mainly concern $\mathbf{Q}$ with codimension $n\geq 2$ of this paper. However, for the convenience of induction and calculation in Section \ref{section4}, we still assume that $n\geq 1$ in Case (1) above. When $n=1$, the only example is $\mathbf{Q}(\xi)=(\xi_1^2)$ by noting $\max_j w_j\not=1$. Then classical restriction estimates for the parabola (see \cite{Zyg74}) gave (\ref{main thm eq2}). 

To illustrate the motivation behind our iterative approach for higher codimensional quadratic surfaces, we  first recall the treatment for classical parabolic case, that is, $\mathbf{Q}(\xi)=(\xi_1^2+\cdots+\xi_d^2)$ with $n=1$. Divide $[0,1]^d$ into squares $\tau$ of $K^{-1}$-scale, and  write $f=\sum_\tau f_\tau =\sum_\tau f\chi_\tau$. We say $\tau_1$ and $\tau_2$ are ``transverse'' if 
\begin{equation}\label{s1a0}
	{\rm dist}(\xi,\eta)\geq 10K^{-1}, \quad \quad \forall ~ \xi\in\tau_1,~ \eta\in \tau_2. 
\end{equation}
Its geometric meaning is that the normal directions of $S_{\mathbf{Q}}$ at the points $\xi$ and $\eta$ are separated. Based on whether or not such a transversality condition is met, we have the following dichotomy:
\begin{align}
    \lvert E^{\mathbf{Q}}f(x)\rvert \lesssim K^{O(1)} \max_{\substack{\tau_1,\tau_2: \\ \text{transverse~}K^{-1}\text{-squares}
    }} | E^{\mathbf{Q}}f_{\tau_1}(x) E^{\mathbf{Q}}f_{\tau_2}(x)|^{\frac{1}{2}}  +\max_{\tau}\Big| E^{\mathbf{Q}}f_{\tau}(x)\Big|.   \label{s1a1}
\end{align}
For the last term in (\ref{s1a1}), note that only a single $\tau$ takes effect, so isotropic rescaling ($\mathbf{Q}$ is quadratically homogeneous) and induction on scale come into play effectively. This fact indicates that the restriction estimate is essentially equivalent to the bilinear restriction estimate for the parabolic case. The research on bilinear restriction estimates is abundant, and the reader can refer to \cite{Gut16,Tao03}.

Next we consider higher codimensional quadratic surfaces with $n\geq 2$. A natural generalization of the transversality condition is that normal planes of $S_{\mathbf{Q}}$ at two points are ``separated,'' which can be captured via $J_{\mathbf{Q}}(\xi;i_1,...,i_{d-n})$. Let us make the geometric intuition precise. Recall that the normal space of $S_\mathbf{Q}$ at $\xi$ is $N_\xi \coloneqq \{(-\nabla\mathbf{Q}(\xi)^T h, h): h\in\mathbb{R}^n \}$. When $d=n$, we have $J_\mathbf{Q}(\xi) = |\nabla\mathbf{Q}(\xi)| = 0$ if and only if there exists $h \neq 0$ such that $\nabla\mathbf{Q}(\xi)^T h = 0$, i.e., $N_\xi\cap N_0 \neq \{0\}$. By passing this relation to other points via the translation-dilation invariant (TDI) property of $S_\mathbf{Q}$, we may call $N_{\xi_1}$ and $N_{\xi_2}$ ``separated'' if $N_{\xi_1} \cap N_{\xi_2} = \{0\}$. Thus, $J_\mathbf{Q}(\xi) \equiv 0$ means that $N_{\xi_1} \cap N_{\xi_2} \neq \{0\}$ whenever $\xi_1 \neq \xi_2$. When $d>n$, things become more complicated because of the parameters $i_1,...,i_{d-n}$ (assume they are all distinct without loss of generality). Via similar considerations, if we define $V_\xi \coloneqq \langle (e_{i_1},0), ..., (e_{i_{d-n}},0) \rangle \oplus N_\xi$ ($\dim V_\xi = d$), then $J_\mathbf{Q}(\xi; i_{1}, ..., i_{d-n}) = 0$ if and only if 
$N_\xi \cap V_0 \neq \{0\}$.
Still by the TDI property of $S_\mathbf{Q}$, we may call $N_{\xi_1}$ and $N_{\xi_2}$ ``separated'' (under $i_1,...,i_{d-n}$) if $V_{\xi_1} \cap N_{\xi_2} = \{0\}$. Thus, $J_\mathbf{Q}(\xi;i_1,...,i_{d-n}) \equiv 0$ means that $V_{\xi_1} \cap N_{\xi_2} \neq \{0\}$ whenever $\xi_1 \neq \xi_2$.

To summarize, for $N_{\xi_1}$ and $N_{\xi_2}$ to be ``separated,'' if $d=n$, then we only need $N_{\xi_1} \cap  N_{\xi_2} = \{0\}$ ($N_{\xi_1} \oplus N_{\xi_2} = \R^{d+n}$); while if $d>n$, then we need the stronger $V_{\xi_1} \cap N_{\xi_2} = \{0\}$ ($V_{\xi_1} \oplus N_{\xi_2} = \R^{d+n}$) under some $i_1,...,i_{d-n}$. On the contrary, if $J_\mathbf{Q}(\xi;i_1,...,i_{d-n})\equiv0$ for any choice of $i_1,...,i_{d-n}$, all pairs of $N_{\xi_1}$ and $N_{\xi_2}$ are ``not separated,'' that is, $V_{\xi_1} \cap N_{\xi_2} \neq \{0\}$ under any $i_1,...,i_{d-n}$. This seemingly unusual geometric interpretation actually makes sense: If $d>n$, then $\{W_1\in \mathrm{Gr}(n,d+n): \dim(W_1 \cap W_2) \geq 1\}$ ($W_2\in \mathrm{Gr}(n,d+n)$ fixed) forms a Schubert variety of higher codimension, which means that its complement might be too ``loose'' to yield useful constraints; on the other hand, extending $N_{\xi_1}$ to some $V_{\xi_1}$ reduces the codimension of the variety to $1$, and so effectively captures curvature properties. The discussions above suggest that $J_\mathbf{Q}$ can be viewed as a measure of certain “generalized minor angles” between normal planes.

In general, one cannot expect $J_\mathbf{Q}$ to be everywhere nonvanishing, and its zero set actually encodes the bilinear transversality property. For example, suppose that $\mathbf{Q}(\xi)=(\xi_1^2,\xi_2^2+\xi_1\xi_3)$ with $d=3$ and $n=2$. Take $i_1=2$, and we get
  $$   |J_{\mathbf{Q}}(\xi;2)| \sim |\xi_1|^2.   $$
  Then we could define a transversality condition that the first coordinates of  two $K^{-1}$-cubes $\tau_1$ and $\tau_2$ are separated. It is worth noting that, since $i_1,...,i_{d-n}$ are all free parameters, we may have multiple ways to define transversality. Still considering the above example, by taking $i_1=3$, we get 
  $$   |J_{\mathbf{Q}}(\xi;3)| \sim |\xi_1||\xi_2|.   $$
  Thus, the separation of the first two coordinates of  two $K^{-1}$-cubes $\tau_1$ and $\tau_2$ can also serve as a suitable transversality condition. However, no matter which transversality condition we consider, it always leads to the following inequality: 
    \begin{align}
  	\lvert E^{\mathbf{Q}}f(x)\rvert \lesssim K^{O(1)} \max_{\substack{\tau_1,\tau_2: \\ \text{transverse~}K^{-1}\text{-squares}
  	}} | E^{\mathbf{Q}}f_{\tau_1}(x) E^{\mathbf{Q}}f_{\tau_2}(x)|^{\frac{1}{2}}  +\max_{W}\Big| \sum_{\tau \subset W}E^{\mathbf{Q}}f_{\tau}(x)\Big|, \label{s1a2}
  \end{align}
  where $\tau_1$ and $\tau_2$ are $K^{-1}$-squares, and $W$ denotes a rectangular box of size $K^{-1}\times 1\times 1$ (and another rectangular box of size $1\times K^{-1}\times 1$ if we consider the second transversality condition) parallel to the coordinate axes. For the bilinear part in (\ref{s1a2}), since $\mathbf{Q}$ is degenerate, classical $L^4$ method is enough to give a nice bound. However, proving the equivalence between the linear restriction estimation and the bilinear restriction estimation for such $\mathbf{Q}$ is no longer trivial. This is because, in the last term in (\ref{s1a2}), it may contain many $\tau$. A standard strategy is to apply decoupling inequalities and proceed as in the parabolic case. Unfortunately, this cannot lead to a sharp bound due to the loss from decoupling. However, we cannot directly use anisotropic rescaling to turn $W$ into $[0,1]^3$ either, since the relevant surface $S_{\mathbf{Q}}$ will change, and induction on scale cannot work. The lack of rescaling invariance is the main difficulty in higher codimensional situations that leads to the loss in restriction estimates.
  
   To overcome this problem, Guo--Oh \cite{GO22} made interesting attempts. For $\mathbf{Q}(\xi)=(\xi_1^2,\xi_2^2+\xi_1\xi_3)$, they excavated more structural information under a nested induction scheme, and derived the sharp restriction estimate. In this paper, we study this type of degenerate quadratic forms from different perspectives. Suppose for the moment that for an $n$-tuple $\mathbf{Q}$ of quadratic forms on $\mathbb{R}^d$, there exist $i_1,...,i_{d-n}$ such that 
   \begin{equation}\label{s1aadd}
   	|J_{\mathbf{Q}}(\xi;i_1,...,i_{d-n})|\sim \prod_{j=1}^t |\xi_{j}|^{s_{j}},
   \end{equation}
   with $s_j\in\N^+$ {\rm(}$1\leq j\leq t${\rm)} and $\sum_{j=1}^t s_j=n.$ In fact, in Section~\ref{section2}, we will prove that all scenarios considered in Theorem \ref{main thm} satisfy this property (after renumbering the variables if necessary). We define $\tau_1$ and $\tau_2$ to be ``transverse'' if the first $t$ coordinates of  any two points $\xi\in \tau_1$ and $\eta\in \tau_2$ are all $10K^{-1}$-separated. Based on this definition, we can use the broad--narrow analysis to obtain (\ref{s1a2}), where $W$ denotes a rectangular box of size $K^{-1}\times 1\times \cdots\times 1$, and its $K^{-1}$ edges are parallel to $e_j$ for some $1\leq j\leq t$. Furthermore, we can iterate (\ref{s1a2}) down to $R^{-1/2}$-scale, then
   \begin{align}
   	\lvert E^{\mathbf{Q}}f(x)\rvert \lesssim K^{O(\frac{K}{\log K})} \sup_{1<\mu_1,...,\mu_t\leq R^{1/2}} \max_{\substack{\tau_1,\tau_2 \\ \text{~are~transverse}}} | E^{\mathbf{Q}}f_{\tau_1}(x) E^{\mathbf{Q}}f_{\tau_2}(x)|^{\frac{1}{2}}  +K^{O(\frac{K}{\log K})}\max_{L}\Big| E^{\mathbf{Q}}f_{L}(x)\Big|,   \label{s1a3}
   \end{align}
   for some $1\leq t\leq d$ (depending on $\mathbf{Q}$), where $\tau_1,\tau_2$ denote two parallel $\mu_1^{-1}\times \cdots \times \mu_t^{-1}\times 1\times \cdots\times 1$ rectangular boxes that are $10\mu_j^{-1}$-separated along the direction of their $\mu_{j}^{-1}$ sides for each $j$, and $L$ is contained in a rectangular box of size $R^{-1/2}\times 1\times \cdots \times 1$, and its $R^{-1/2}$ edges are parallel to $e_j$ for some $1\leq j\leq t$.

   Note that $f_{L}$ in (\ref{s1a3}) is a genuinely lower dimensional object, since it is supported at the (smallest possible) scale $R^{-1/2}$ rather than $K^{-1}$. This change gives us more freedom to choose other techniques, such as the locally constant property, adaptive decoupling and adaptive restriction estimates, not just induction on scale.  By applying these tools on a case-by-case basis, we can get optimal bounds for the linear term in (\ref{s1a3}). 
   
   On the other hand, the difficulty of rescaling mentioned above has now shifted to the first bilinear term in (\ref{s1a3}). There are $t$ uncontrollable transversality parameters $\mu_1,...,\mu_t$, and if we use anisotropic rescaling to change all of them to $\sim 1$, then there is no guarantee that $S_{\mathbf{Q}}$ is preserved, and so classical bilinear restriction estimates for such $\mathbf{Q}$ are not applicable. We also point out that for the parabolic case, such argument can work and lead to the same bound derived via (\ref{s1a1}), see \cite[section 4]{demeter2020}.
   
   To overcome this obstacle, we will build a stronger bilinear restriction estimate (see Theorem~\ref{s3thm1}) that is uniform over all transversality parameters. A notable fact is that the endpoints of such estimates rely on $\max_j s_j$, where $s_j$ is defined in (\ref{s1aadd}) and may depend on our selection of $i_1,...,i_{d-n}$. We take $\mathbf{Q}(\xi)=(\xi_1^2,\xi_2^2+\xi_1\xi_3)$ as an example. From the above computations, we know $\max_j s_j=2$ and $1$ when $i_1=2$ and $3$, respectively. The latter can provide us with a wider range and help us obtain the desired bound, while the former cannot. Therefore, to obtain acceptable bounds for bilinear terms, we may need to choose $i_1,...,i_{d-n}$ very carefully to lower $\max_j s_j$ as much as possible. Fortunately, this can be done by suitably adjusting the parameters $i_1,...,i_{d-n}$ (Theorem~\ref{thm:Jacobian}) while retaining all the desired properties of $J_\mathbf{Q}$. So, the remaining issue is why we could have the form (\ref{s1aadd}), since Theorem~\ref{main thm} itself only assumes ``$J_\mathbf{Q}\not\equiv0$''. It might be a surprise that such a linear factorization can be guaranteed merely by the nonvanishing assumption (Theorem~\ref{thm:Jacobian}). 
   
   The proof of Theorem~\ref{thm:Jacobian} relies on tools and techniques from both algebra and graph theory. In fact, we could even get an almost complete characterization of those monomial-type $\mathbf{Q}$ with $J_\mathbf{Q}\equiv0$: Roughly speaking, if $J_\mathbf{Q}(\xi;i_1,...,i_{d-n}) \equiv 0$ for all choices of $i_1,...,i_{d-n}$, then there exists a ``cycle'' structure with even length. More precisely, 
   $\mathbf{Q}(\xi) = (\xi_1\xi_2, \xi_2\xi_3, ..., \xi_{2m-1}\xi_{2m}, \xi_{2m}\xi_1)$ with $m\geq2$ has $J_\mathbf{Q} \equiv 0$, and our discussions indicate that this is in some sense the only nontrivial way that such degeneracy can happen when $d \geq n$. See Section~\ref{section2} for more comments on our proof strategy.

\vskip0.3cm

\noindent \textbf{Outline of the paper.}  In Section \ref{section2}, we prove structural properties of $J_{\mathbf{Q}}$ via algebra and graph theory. In Section \ref{section3}, we first give a detailed exposition of the broad--narrow analysis and iteration for general higher codimensional quadratic surfaces, and then establish general bilinear restriction estimates independent of the transversality parameters. In Section \ref{section4}, we present the main proof of Theorem \ref{main thm}. In Section \ref{section5}, we provide additional remarks and applications of our method. 

\vskip0.3cm

\noindent \textbf{Notation.} We will use $\# X$ to denote the cardinality of a finite set $X$, and use $|X|$ to denote the Lebesgue measure of a measurable set $X$. We will use $B_R^m$ to represent a ball with radius $R$ in $\mathbb{R}^m$, and we abbreviate $B_R^m$ as $B_R$ if $m$ is clear from the context. We will write $A\lesssim_\epsilon B$ to mean that there exists a constant $C$ depending on $\epsilon$ such that $A\leq CB$. Moreover, $A \sim B$ means $A \lesssim B$ and $A\gtrsim B$. Let $e(b):=e^{2\pi i b}$ for each $b\in \mathbb{R}$.  
We will use $\on{supp}f$ to denote the support of a function $f$, and use $\chi_E$ to denote the characteristic function of a set $E$. Define $f_E:=f\chi_E$ for any function $f$ and any set $E$. Let $p'$ denote the dual exponent of $p$, that is, $1/p+1/p'=1$.
Let $I_{n\times n}$ be the standard identity matrix of order $n$ for any $n\in\mathbb{N}^+$. Let $\N$ denote the set of natural numbers including zero. Whenever we say ``on $\mathbb{R}^D$'' for some $D$, it is implicitly assumed that the tuple in question depends on exactly $D$ variables and does not miss any of them.

\section{Properties of the generalized Jacobian}\label{section2}

The main result of this section is as follows:
\begin{theorem}\label{thm:Jacobian}
Let $d \geq n$, and $\mathbf{Q}(\xi) = (\xi_{\lambda_{2m-1}} \xi_{\lambda_{2m}})_{m=1}^n$ be an $n$-tuple of quadratic monomials on $\mathbb{R}^d$. Suppose that $J_{\mathbf{Q}}(\xi;i_1,...,i_{d-n}) \not\equiv 0$ for some $i_1,...,i_{d-n}$, then we must have 
\begin{equation}\label{eq:linear_factor}
   	J_{\mathbf{Q}}(\xi;i_1,...,i_{d-n}) = c\cdot \prod_{j=1}^d\xi_{j}^{s_j}
\end{equation}
for some $c\neq0$, where $s_j\in\N$ {\rm(}$1\leq j\leq d${\rm)} and $\sum_{j=1}^d s_{j}=n$. 

Moreover, let $w_j$ {\rm(}$1\leq j\leq d${\rm)} be the number of occurrences of $\xi_j$ in $\{ \xi_{\lambda_j}\}_{j=1}^{2n}$, then $\max_j s_j \leq \max_j w_j$. If $\max_j w_j \neq 1$, then there exist $i_1',...,i_{d-n}'$ such that $J_{\mathbf{Q}}(\xi;i_1',...,i_{d-n}') \not\equiv 0$ has the form ${\rm(\ref{eq:linear_factor})}$ with $s_j$ replaced by $s_j'$, and 
\begin{align}\label{eq:degree_ineq}
    \max_j s_j' \leq \max_j w_j - 1.
\end{align}
\end{theorem}

Theorem~\ref{thm:Jacobian} is vital for the proof of Theorem~\ref{main thm} in three aspects. First, it basically tells us that $J_\mathbf{Q}$ enjoys a linear factorization as long as it is not identically zero, which induces the transversality condition and the associated uniform bilinear estimates (see the argument of Section \ref{section3}). Second, the by-product $\max_j s_j' \leq \max_j w_j - 1$ allows us to show that bilinear estimates are always acceptable, thus reducing everything to the analysis of linear terms (see the argument of Section \ref{section4}). One caveat is that if we are not allowed to adjust the original parameters $i_1,...,i_{d-n}$, then this may not hold true in general. Finally, its proof and the ideas involved are much stronger than the theorem itself, in the sense that we could get an almost complete characterization of those $\mathbf{Q}$ with $J_{\mathbf{Q}}(\xi;i_1,...,i_{d-n}) \not\equiv 0$ for some $i_1,...,i_{d-n}$. This makes it possible to obtain other useful properties that we will need later on, such as Theorem~\ref{thm:inherit}.

The proof of Theorem~\ref{thm:Jacobian} relies on tools and techniques from both algebra and graph theory. Such an approach does not seem to appear in any previous research of higher codimensional restriction estimates, and may be of interest in itself. The rough idea is to construct a graph encoding the algebraic structure of $\mathbf{Q}$ and reinterpret everything in the language of graph theory. For example, nonzero terms in the permutation expansion of $J_{\mathbf{Q}}(\xi;i_1,...,i_{d-n})$ correspond to ``admissible'' digraph structures, so the analysis becomes much easier and more intuitive. Now we provide the details. Our notation and terminology closely follow those of \cite{diestel25}, especially sections~1.1--1.5 therein.

\begin{proof}[Proof of Theorem~\ref{thm:Jacobian}]
Without loss of generality, we may always assume that the quadratic monomials in $\mathbf{Q}$ are all distinct. Otherwise, it is easy to see that two rows of $\nabla\mathbf{Q}$ coincide, and thus $J_{\mathbf{Q}}(\xi;i_1,...,i_{d-n}) \equiv 0$ for any $i_1,...,i_{d-n}$, a contradiction. Similarly, we may assume that $\mathbf{Q}$ does not miss any variable.

Construct a directed graph $G=(V,E)$ associated with  $\mathbf{Q}$ with vertices $V=\{\xi_j\}_{j=1}^d$, edges $E= \{\xi_{\lambda_{2m-1}}\xi_{\lambda_{2m}}\}_{m=1}^n$, and two direction maps init: $E \rightarrow V$ and ter: $E \rightarrow V$ (see \cite[section~1.10]{diestel25}). By the distinctness, $E$ faithfully encodes all terms in $\mathbf{Q}$. For $e=\xi_{\lambda_{2m-1}}\xi_{\lambda_{2m}} \in E$, we have ${\rm init}(e) = {\rm ter}(e)=\xi_{\lambda_{2m}}$ (fixed loop) if $e$ is squared, and $\{{\rm init}(e), {\rm ter}(e)\} = \{\xi_{\lambda_{2m-1}}, \xi_{\lambda_{2m}}\}$ (undetermined edge) if $e$ is mixed. For later convenience, we define: 
\begin{itemize}
    \item $V_0 := \{\xi_j\in V: \xi_j^2 \in E\}$ (vertices of loops).
    \item $W_0 := \{\xi_j\in V: j \in \{i_1,...,i_{d-n}\} \}$ (tagged vertices corresponding to the parameters).
    \item ${\rm init}^{-1}(\xi_j) := \{ e\in E: {\rm init}(e) = \xi_j\}$ (outgoing edges of a vertex $\xi_j$).
    \item ${\rm ter}^{-1}(\xi_j) := \{ e\in E: {\rm ter}(e) = \xi_j\}$ (incoming edges of a vertex $\xi_j$).
\end{itemize}
We call a pair (init, ter) \textit{admissible} if $\# {\rm init}^{-1}(\xi_j) = 1$ for each $\xi_j\notin W_0$, and $\# {\rm init}^{-1}(\xi_j)=0$ for each $\xi_j \in W_0$. 

Since any nonzero term in the permutation expansion of $J_{\mathbf{Q}}(\xi;i_1,...,i_{d-n})$ picks exactly one entry in each row, it naturally induces a pair (init,\,ter) as follows. For $e=\xi_{\lambda_{2m-1}}\xi_{\lambda_{2m}} \in E$, if $\xi_{\lambda_{2m}}$ is picked in the $m$$j$th row, then let $({\rm init}(e), {\rm ter}(e))=(\xi_{\lambda_{2m-1}}, \xi_{\lambda_{2m}})$; while if $\xi_{\lambda_{2m-1}}$ is picked in the $m$$j$th row, then let $({\rm init}(e), {\rm ter}(e))=(\xi_{\lambda_{2m}}, \xi_{\lambda_{2m-1}})$. In other words, $({\rm init}(e), {\rm ter}(e))=(\xi_{\lambda_{2m-1}}, \xi_{\lambda_{2m}})$ means that we ``pick $\xi_{\lambda_{2m}}$ in the $\lambda_{2m-1}$$j$th column,'' and vice versa. Since the nonzero term also picks exactly one entry in each column, it is easy to check that its induced (init, ter) must be admissible. Conversely, any admissible pair (init,\,ter) uniquely determines a nonzero term in the permutation expansion: For each $\xi_j\notin W_0$, since $\xi_j = {\rm init}(e)$ for a unique $e\in E$ by admissibility, we pick ${\rm ter}(e)$ in the $j$th column. Therefore, we get a one-to-one correspondence between nonzero terms and admissible pairs. In particular, $J_{\mathbf{Q}}(\xi;i_1,...,i_{d-n}) \not\equiv 0$ means that there exists at least one admissible pair (init,\,ter).

For convenience, we will digress to discuss the relationship between $\max_j s_j$ ($\max_j s_j'$) and $\max_j w_j$ whenever appropriate during our proof of (\ref{eq:linear_factor}), and we are assuming the truth of (\ref{eq:linear_factor}) for these discussions. However, one can see that the proof of (\ref{eq:linear_factor}) itself is self-contained, so there will be no circular reasoning. 

Note that $w_j = \# {\rm init}^{-1}(\xi_j) + \# {\rm ter}^{-1}(\xi_j)$ (the degree of the vertex $\xi_j$ in the digraph $G$) and $s_j \leq \max_{\rm (init,\,ter)}\# {\rm ter}^{-1}(\xi_j)$, where ``$\max$'' is taken over the set of all admissible pairs (init,\,ter). So, we must have $s_j \leq w_j$ for each $j$, which means $\max_j s_j \leq \max_j w_j$, as desired. Also note that since $\# {\rm init}^{-1}(\xi_j) = 1$ for each $\xi_j\in V_0$ by definition, we must have $V_0 \cap W_0 = \varnothing$.

Decompose $G$ into its components (i.e., maximal connected subgraphs) $G_k=(V_k,E_k)$, $1\leq k \leq K$. We claim $\# V_k \cap (V_0 \sqcup W_0) \leq 1$. Otherwise, by the connectivity of $G_k$, we can always find a path $P=\xi_{j_0}\xi_{j_1}\cdots\xi_{j_L}$ ($L\in\N^+$) between $\xi_{j_0} \neq \xi_{j_L} \in V_0 \sqcup W_0$ in $G_k$. Without loss of generality, we may take $P$ to be the shortest such path, so that $\xi_{j_l} \notin V_0 \sqcup W_0$ for each $1\leq l \leq L-1$. Fix any admissible pair (init, ter). Note that $\xi_{j_0} \in V_0 \sqcup W_0$ implies ${\rm ter}(\xi_{j_0}\xi_{j_1}) = \xi_{j_0}$ and ${\rm init}(\xi_{j_0}\xi_{j_1}) = \xi_{j_1}$ by definition, which in turn implies ${\rm ter}(\xi_{j_1}\xi_{j_2}) = \xi_{j_1}$ and ${\rm init}(\xi_{j_1}\xi_{j_2}) = \xi_{j_2}$. Proceeding inductively, we eventually have ${\rm ter}(\xi_{j_{L-1}}\xi_{j_L}) = \xi_{j_{L-1}}$ and ${\rm init}(\xi_{j_{L-1}}\xi_{j_L}) = \xi_{j_L}$. However, $\xi_{j_L} \in V_0 \sqcup W_0$ implies ${\rm ter}(\xi_{j_{L-1}}\xi_{j_L}) = \xi_{j_{L}}$ and ${\rm init}(\xi_{j_{L-1}}\xi_{j_L}) = \xi_{j_{L-1}}$, a contradiction. So, the claim holds true.

To summarize, we can divide the components $\{G_k\}_{k=1}^K$ into three types: 
\begin{enumerate}[label={\rm (\roman*)}]
    \item $V_k \cap (V_0 \sqcup W_0) = \varnothing$ (no loop or tagged vertex);
    \item $V_k \cap W_0 = \varnothing$ and $\# V_k \cap V_0 = 1$ (no tagged vertex, exactly one loop);
    \item $V_k \cap V_0 = \varnothing$ and $\# V_k \cap W_0 = 1$ (no loop, exactly one tagged vertex).
\end{enumerate}
In particular, this implies $K \geq \# W_0 = d-n$.

Since $V=\sqcup_k V_k$ and $E=\sqcup_k E_k$, the Jacobian matrix in question enjoys a special block form, where each block corresponds to exactly one $G_k$. Now we make this intuition precise. We claim the following. 
\begin{itemize}
    \item If $G_k$ has Type ${\rm (i)}$ or ${\rm (ii)}$, then $\# V_k = \# E_k$.
\end{itemize}
To see this, first note that there is a block with columns indexed by $V_k$ and rows indexed by $E_k$. Assume, to the contrary, that there exists $G_k$ of Type ${\rm (i)}$ with $\# V_k \neq \# E_k$. If $\# V_k > \# E_k$, then by a Laplace expansion along the columns corresponding to $V_k$, we get $J_{\mathbf{Q}}(\xi;i_1,...,i_{d-n}) \equiv 0$, a contradiction. On the other hand, if $\# V_k < \# E_k$, then by a Laplace expansion along the rows corresponding to $E_k$, we still get $J_{\mathbf{Q}}(\xi;i_1,...,i_{d-n}) \equiv 0$, a contradiction. So, the claim holds true. Similar arguments can also be used to show the following. 
\begin{itemize}
    \item If $G_k$ has Type ${\rm (iii)}$, then $\# V_k = \# E_k +1$.
\end{itemize}
Here ``$+1$'' comes from adding the $(n+t)$$j$th row $e_{i_t}$ of the Jacobian matrix to the block if $V_k \cap W_0 = \{\xi_{i_t}\}$ ($1 \leq t \leq d-n$). 
Therefore, the square block $M_k$ associated to $G_k$ has size $\# E_k$ by $\# V_k$ for $G_k$ of Type ${\rm (i)}$ or ${\rm (ii)}$, and size $(\# E_k + 1)$ by $\# V_k$ for $G_k$ of Type ${\rm (iii)}$. Also, since $J_{\mathbf{Q}}(\xi;i_1,...,i_{d-n}) = |M_1|\cdots|M_K| \not\equiv 0$, we have $|M_k| \not\equiv 0$ for each $1\leq k \leq K$. Thus, the original problem is reduced to simpler subproblems of the same type. More precisely, it suffices to prove (\ref{eq:linear_factor}) and (\ref{eq:degree_ineq}) for connected $G=(V,E)$ (constructed as before) in two ``special cases''.
\begin{enumerate}[label={\rm (\Roman*)}]
    \item $d=n$ (i.e., $\# V = \# E$), $G$ has no loop or exactly one loop.
    \begin{itemize}
        \item 
        Suppose that $|\nabla \mathbf{Q}| \not\equiv 0$, then $|\nabla \mathbf{Q}|$ has the form (\ref{eq:linear_factor}).
        \item 
        Moreover, $\max_j w_j \neq 1$, and $\max_j s_j \leq \max_j w_j - 1$.
    \end{itemize}
    
    \item $d=n+1$ (i.e., $\# V = \# E + 1$), G has no loop.
    \begin{itemize}
        \item 
        For any $i$, we always have $J_\mathbf{Q}(\xi;i) \not\equiv 0$, and $J_\mathbf{Q}(\xi;i)$ has the form (\ref{eq:linear_factor}).
        \item 
        Moreover, if $\max_j w_j \neq 1$, then there exists $i'$ such that $J_\mathbf{Q}(\xi;i')\not\equiv0$ has the form (\ref{eq:linear_factor}) with $s_j$ replaced by $s_j'$, and $\max_j s_j' \leq \max_j w_j - 1$.
    \end{itemize}
\end{enumerate}
This is because (\ref{eq:linear_factor}) and (\ref{eq:degree_ineq}) are ``inherited upward'': As long as we have established their natural analogs for each $G_k$, which corresponds to one of the above cases, they will also hold for $G$. The details of this reduction are simple, and so are left to the reader. We only point out that the second conclusion in Case (I) and the first conclusion in Case (II) are actually stronger than those in Theorem~\ref{thm:Jacobian} respectively. Now we handle Case (I) and Case (II), respectively. 

For Case (I), since $|\nabla\mathbf{Q}|\not\equiv0$, there exists at least one admissible pair (init,\,ter). And for any such pair, we must have $\# {\rm init}^{-1}(\xi_j) = 1$ for each $\xi_j \in V$. Recall that $w_j = \# {\rm init}^{-1}(\xi_j) + \# {\rm ter}^{-1}(\xi_j)$ and $s_j \leq \max_{\rm (init,\,ter)}\# {\rm ter}^{-1}(\xi_j)$, so we must have $s_j \leq w_j - 1$ for each $j$, which means $\max_j s_j \leq \max_j w_j - 1$, as desired. In particular, this implies $\max_j w_j \neq 1$ as $\max_j s_j \geq 1$.

When $G$ has no loop, by \cite[Corollary~1.5.2]{diestel25}, $G$ cannot be a tree, that is, $G$ contains a cycle $C=\xi_{j_1}\xi_{j_2}\cdots\xi_{j_L}\xi_{j_1}$ ($L \geq 3$). Note that $G-\xi_{j_L}\xi_{j_1}$ is still connected and has $\# V = \# (E-\xi_{j_L}\xi_{j_1}) + 1$, and so must be a tree by \cite[Corollary~1.5.2]{diestel25} again. By \cite[Theorem 1.5.1 (ii)]{diestel25}, $\xi_{j_1}\xi_{j_2}\cdots\xi_{j_L}$ is the unique path in $G-\xi_{j_L}\xi_{j_1}$, and so $C$ is the only cycle in $G$. Besides, it is easy to check that either ${\rm init}(\xi_{j_l}\xi_{j_{l+1}}) = \xi_{j_l}$ (here $\xi_{j_{L+1}}\equiv \xi_{j_{1}}$) for each $1\leq l \leq L$, or ${\rm init}(\xi_{j_l}\xi_{j_{l+1}}) = \xi_{j_{l+1}}$ for each $1\leq l \leq L$.

If $L=d$, then $G=C$, and $\mathbf{Q} = (\xi_{j_1}\xi_{j_2}, \xi_{j_2}\xi_{j_3}, ..., \xi_{j_d}\xi_{j_1}) \equiv (\xi_1\xi_2, \xi_2\xi_3, ..., \xi_d\xi_1)$ up to a harmless relabeling of the vertices. Thus, by direct computation, $|\nabla\mathbf{Q}| = \pm 2\xi_1\xi_2\cdots\xi_d$ for odd $d$, and $|\nabla\mathbf{Q}| \equiv 0$ for even $d$. Since $|\nabla\mathbf{Q}| \not\equiv0$, we must have $d$ is odd and $|\nabla\mathbf{Q}|$ has the form (\ref{eq:linear_factor}), as desired. If $L < d$, then each component of $G - C$ must be a tree intersecting $C$ at exactly one vertex, which we define to be the root of the tree. For each $1\leq l \leq L$, let $T_l$ be the tree rooted at $\xi_{j_l}$. Note that different trees are disjoint, and we allow the trivial case $T_l = \{\xi_{j_l}\}$. Let $\Lambda$ be the set of $l$ such that $T_l$ is nontrivial. For each $l\in\Lambda$, let $\xi^{(l)}$ be the vector composed of all vertices of $T_l$, with $\xi_{j_l}$ being its first coordinate, and $\mathbf{Q}_l(\xi^{(l)})$ be the quadratic form composed of all edges of $T_l$. Let $\mathbf{Q}_0 := (\xi_{j_1}\xi_{j_2}, \xi_{j_2}\xi_{j_3}, ..., \xi_{j_L}\xi_{j_1})$. The key observation is that $$|\nabla \mathbf{Q}| = \pm |\nabla\mathbf{Q}_0|\prod_{l\in\Lambda} J_{\mathbf{Q}_l}(\xi^{(l)}; 1).$$
Since $|\nabla \mathbf{Q}|\not\equiv0$, we must have $|\nabla \mathbf{Q}_0|\not\equiv0$. Via similar arguments as in the $L=d$ case, we must have $L$ is odd and $|\nabla \mathbf{Q}_0| = \pm 2\xi_{j_1}\xi_{j_2}\cdots\xi_{j_L}$. On the other hand, each $T_l$ exactly satisfies the assumptions in Case (II), which means $J_{\mathbf{Q}_l}(\xi^{(l)}; 1)$ has the form (\ref{eq:linear_factor}) for each $l\in\Lambda$. Therefore, granting that the conclusions in Case (II) hold, $|\nabla \mathbf{Q}|$ always has the form (\ref{eq:linear_factor}), and so the proof for Case (I) when $G$ has no loop is completed.

When $G$ has exactly one loop $\xi_j^2$, let $G':= G - \xi_j^2$, then $G'$ is still connected and has $\# V = \# (E - \xi_j^2) + 1$. Let $\mathbf{Q}'$ be the quadratic form obtained by deleting $\xi_j^2$ from $\mathbf{Q}$. The key observation is that $|\nabla\mathbf{Q}| = \pm 2\xi_j J_{\mathbf{Q}'}(\xi;j)$. Note also that $G'$ has no loop, so it exactly satisfies the assumptions in Case (II), which means $J_{\mathbf{Q}'}(\xi;j)$ has the form (\ref{eq:linear_factor}). Therefore, granting that the conclusions in Case (II) hold, $|\nabla \mathbf{Q}|$ always has the form (\ref{eq:linear_factor}). So, the proof for Case (I) when $G$ has exactly one loop, thus the whole case (I), is completed. Now everything comes down to handling Case (II).

For Case (II), by \cite[Corollary~1.5.2]{diestel25}, $G$ must be a tree. For any $i$, let $\xi_i$ be the root of $G$. We claim that there is exactly one nonzero term in the permutation expansion of $J_\mathbf{Q}(\xi;i)$, and thus the first conclusion immediately holds true. To prove the claim, it suffices to show that there is exactly one admissible pair (init,\,ter). By the definition of admissibility, we need $\# {\rm init}^{-1}(\xi_j) = 1$ for $j\neq i$, and $\# {\rm init}^{-1}(\xi_i)=0$. Thanks to \cite[Theorem 1.5.1]{diestel25}, we can properly define the tree-order associated with $G$ and $\xi_i$, as in \cite[p. 15]{diestel25}. By induction on height of vertices, it is easy to check that the only possible construction of (init,\,ter) is: For any two adjacent vertices $\xi_j$ and $\xi_{j'}$ in $G$, define $({\rm init}(\xi_j\xi_{j'}), {\rm ter}(\xi_j\xi_{j'})) = (\xi_j,\xi_{j'})$ if and only if $\xi_j > \xi_{j'}$ in the tree-order. So, the claim holds true.

For the second claim in Case (II), simply take $i'$ so that $\xi_{i'}$ is a leaf of $G$. Since $G \neq \varnothing$, there exists at least one such $i$. By the first conclusion, we must have $J_\mathbf{Q}(\xi;i') \not\equiv 0$, and $J_\mathbf{Q}(\xi;i')$ has the form (\ref{eq:linear_factor}) with $s_j$ replaced by $s_j'$. Moreover, by the proof in the last paragraph, we know that $J_\mathbf{Q}(\xi;i')$ is given by a unique admissible pair (init,\,ter) corresponding to the tree-order associated with $G$ and $\xi_{i'}$. Also, we have $\# {\rm init}^{-1}(\xi_j) = 1$ for $j\neq i'$, and $\# {\rm init}^{-1}(\xi_{i'})=0$. For any $w_{j_0} = \max_j w_j$, since $\max_j w_j \neq 1$ and $w_{i'} = 1$ by definition of leaves, we must have $\xi_{j_0} \neq \xi_{i'}$, which means $\# {\rm init}^{-1}(\xi_{j_0}) = 1$. However, recall that $w_{j_0} = \# {\rm init}^{-1}(\xi_{j_0}) + \# {\rm ter}^{-1}(\xi_{j_0})$ and $s_{j_0}' \leq \# {\rm ter}^{-1}(\xi_{j_0})$, so we must have $s_{j_0}' \leq w_{j_0} - 1$. Therefore, we can conclude that $\max_j s_j' \leq \max_j w_j - 1$, as desired. This completes the proof for Case (II), and thus Theorem~\ref{thm:Jacobian}.
\end{proof}
Before proceeding, we provide several representative examples to help the reader digest the proof of Theorem~\ref{thm:Jacobian}. Figures~\ref{fig:two_cycle}, \ref{fig:cycle_loop}, and $\ref{fig:even_cycle}$ are typical illegal examples. There are also legal examples, and here we only give one of them (Figure~\ref{fig:legal}), which already contains all types of components discussed in the proof.
\begin{figure}[htbp]
  \centering
  
    \hspace*{0.5cm}  
    \tikzset{every picture/.style={line width=0.75pt}} 
    
\begin{tikzpicture}[x=0.75pt,y=0.75pt,yscale=-1,xscale=1]

\draw    (532.1,50.68) -- (511.33,107.16) ;
\draw   (328.71,50.97) -- (328.71,119.2) -- (291.21,88.08) -- cycle ;
\draw    (328.71,119.2) -- (402.65,118.88) ;
\draw    (402.65,118.88) -- (402.65,50.66) ;
\draw   (402.65,118.88) -- (402.65,50.66) -- (440.27,84.97) -- cycle ;

\draw (282.87,91.46) node   [align=left] {$\displaystyle \xi _{1}$};
\draw (322.45,41.36) node   [align=left] {$\displaystyle \xi _{3}$};
\draw (323.94,126.97) node   [align=left] {$\displaystyle \xi _{2}$};
\draw (411.9,126.49) node   [align=left] {$\displaystyle \xi _{4}$};
\draw (411.21,41.86) node   [align=left] {$\displaystyle \xi _{5}$};
\draw (522.92,113.39) node   [align=left] {$\displaystyle \xi _{7}$};
\draw (542.47,59.86) node   [align=left] {$\displaystyle \xi _{8}$};
\draw (449.24,90.25) node   [align=left] {$\displaystyle \xi _{6}$};
\draw (328.71,50.97) node    {$\bullet $};
\draw (291.21,88.08) node    {$\bullet $};
\draw (328.71,119.2) node    {$\bullet $};
\draw (402.65,50.66) node    {$\bullet $};
\draw (440.27,84.97) node    {$\bullet $};
\draw (402.65,118.88) node    {$\bullet $};
\draw (532.1,50.68) node    {$\bullet $};
\draw (511.33,107.16) node    {$\bullet $};
\end{tikzpicture}

  \caption{$\mathbf{Q}(\xi) = (\xi_1\xi_2, \xi_2\xi_3, \xi_3\xi_1, \xi_2\xi_4, 
  \xi_4\xi_5, \xi_5\xi_6,
  \xi_5\xi_4,
  \xi_7\xi_8)$ with $d=n=8$. The left component has two cycles, so $J_\mathbf{Q}(\xi) \equiv 0$.}  
    
  \label{fig:two_cycle}
\end{figure}
\begin{figure}[htbp]
  \centering
  \hspace*{0.5cm}  

\tikzset{every picture/.style={line width=0.75pt}} 

\begin{tikzpicture}[x=0.75pt,y=0.75pt,yscale=-1,xscale=1]

\draw    (336.26,315.64) -- (318.25,362.53) ;
\draw   (318.71,299.23) .. controls (318.71,290.17) and (326.56,282.82) .. (336.26,282.82) .. controls (345.95,282.82) and (353.81,290.17) .. (353.81,299.23) .. controls (353.81,308.29) and (345.95,315.64) .. (336.26,315.64) .. controls (326.56,315.64) and (318.71,308.29) .. (318.71,299.23) -- cycle ;
\draw   (228.1,308.69) -- (318.25,363.78) -- (228.1,363.78) -- cycle ;
\draw    (449.33,304.67) -- (412,362.67) ;
\draw    (487.33,363.67) -- (412,362.67) ;
\draw    (594.33,305.67) -- (557,363.67) ;
\draw    (632.33,364.67) -- (557,363.67) ;

\draw (328.09,370.41) node   [align=left] {$\displaystyle \xi _{2}$};
\draw (346.04,324.79) node   [align=left] {$\displaystyle \xi _{4}$};
\draw (222.38,371.12) node   [align=left] {$\displaystyle \xi _{1}$};
\draw (237.29,299.39) node   [align=left] {$\displaystyle \xi _{3}$};
\draw (228.1,308.69) node    {$\bullet $};
\draw (336.26,315.64) node    {$\bullet $};
\draw (228.1,363.78) node    {$\bullet $};
\draw (318.25,363.78) node    {$\bullet $};
\draw (461.92,295.57) node   [align=left] {$\displaystyle \xi _{7}$};
\draw (402.92,370.57) node   [align=left] {$\displaystyle \xi _{5}$};
\draw (498.38,371.52) node   [align=left] {$\displaystyle \xi _{6}$};
\draw (449.33,304.67) node    {$\bullet $};
\draw (412,362.67) node    {$\bullet $};
\draw (487.33,363.67) node    {$\bullet $};
\draw (606.92,296.57) node   [align=left] {$\displaystyle \xi _{10}$};
\draw (547.92,371.57) node   [align=left] {$\displaystyle \xi _{8}$};
\draw (643.38,372.52) node   [align=left] {$\displaystyle \xi _{9}$};
\draw (594.33,305.67) node    {$\bullet $};
\draw (557,363.67) node    {$\bullet $};
\draw (632.33,364.67) node    {$\bullet $};

\end{tikzpicture}

\caption{$\mathbf{Q}(\xi) = (\xi_1\xi_2, \xi_2\xi_3, \xi_3\xi_1, \xi_2\xi_4, \xi_4^2, \xi_5\xi_6, \xi_5\xi_7, \xi_8\xi_9, \xi_8\xi_{10})$ with $d=10$ and $n=9$. The left component has a cycle and a loop, so $J_\mathbf{Q}(\xi;i) \equiv 0$ for any $i$.}

  \label{fig:cycle_loop}
\end{figure}
\begin{figure}[htbp]
  \centering
  \hspace*{0.5cm}  

\tikzset{every picture/.style={line width=0.75pt}} 

\begin{tikzpicture}[x=0.75pt,y=0.75pt,yscale=-1,xscale=1]

\draw   (264.11,689.46) -- (338.27,689.46) -- (338.27,758.84) -- (264.11,758.84) -- cycle ;
\draw    (459.33,693.49) -- (422,751.49) ;
\draw    (497.33,752.49) -- (422,751.49) ;
\draw    (224.57,721.75) -- (264.11,758.84) ;

\draw (256.82,679.69) node   [align=left] {$\displaystyle \xi _{4}$};
\draw (258.69,766.74) node   [align=left] {$\displaystyle \xi _{1}$};
\draw (348.07,766.25) node   [align=left] {$\displaystyle \xi _{2}$};
\draw (346.2,680.2) node   [align=left] {$\displaystyle \xi _{3}$};
\draw (471.92,684.39) node   [align=left] {$\displaystyle \xi _{8}$};
\draw (412.92,759.39) node   [align=left] {$\displaystyle \xi _{6}$};
\draw (508.38,760.34) node   [align=left] {$\displaystyle \xi _{7}$};
\draw (264.11,689.46) node    {$\bullet $};
\draw (338.27,689.46) node    {$\bullet $};
\draw (264.11,758.84) node    {$\bullet $};
\draw (338.27,758.84) node    {$\bullet $};
\draw (459.33,693.49) node    {$\bullet $};
\draw (422,751.49) node    {$\bullet $};
\draw (497.33,752.49) node    {$\bullet $};
\draw (224.57,721.75) node    {$\bullet $};
\draw (214.35,715.99) node   [align=left] {$\displaystyle \xi _{5}$};

\end{tikzpicture}

\caption{$\mathbf{Q}(\xi) = (\xi_1\xi_2, \xi_2\xi_3, \xi_3\xi_4, \xi_4\xi_1,  
\xi_1\xi_5, \xi_6\xi_7, \xi_6\xi_8)$ with $d=8$ and $n=7$. The left component has a cycle of even length, so $J_\mathbf{Q}(\xi;i) \equiv 0$ for any $i$.}

  \label{fig:even_cycle}
\end{figure}
\begin{figure}[htbp]
  \centering
  \hspace*{0.5cm}  

\tikzset{every picture/.style={line width=0.75pt}} 

\begin{tikzpicture}[x=0.75pt,y=0.75pt,yscale=-1,xscale=1]

\draw   (213,497.16) -- (248,553.16) -- (178,553.16) -- cycle ;
\draw    (295.33,514.49) -- (247.33,553.16) ;
\draw    (388.1,522.68) -- (351.33,550.82) ;
\draw   (367.86,502.92) .. controls (367.86,492) and (376.92,483.16) .. (388.1,483.16) .. controls (399.27,483.16) and (408.33,492) .. (408.33,502.92) .. controls (408.33,513.83) and (399.27,522.68) .. (388.1,522.68) .. controls (376.92,522.68) and (367.86,513.83) .. (367.86,502.92) -- cycle ;
\draw    (506.33,495.67) -- (469,553.67) ;
\draw    (544.33,554.67) -- (469,553.67) ;

\draw (225.92,488.39) node   [align=left] {$\displaystyle \xi _{3}$};
\draw (171.92,561.39) node   [align=left] {$\displaystyle \xi _{1}$};
\draw (258.38,560.34) node   [align=left] {$\displaystyle \xi _{2}$};
\draw (294.47,499.86) node   [align=left] {$\displaystyle \xi _{4}$};
\draw (213,497.16) node    {$\bullet $};
\draw (178,553.16) node    {$\bullet $};
\draw (248,553.16) node    {$\bullet $};
\draw (295.33,514.49) node    {$\bullet $};
\draw (360.71,560.01) node   [align=left] {$\displaystyle \xi _{5}$};
\draw (396.47,532.86) node   [align=left] {$\displaystyle \xi _{6}$};
\draw (388.1,522.68) node    {$\bullet $};
\draw (351.33,550.82) node    {$\bullet $};
\draw (518.92,486.57) node   [align=left] {$\displaystyle \xi _{9}$};
\draw (459.92,561.57) node   [align=left] {$\displaystyle \xi _{7}$};
\draw (555.38,562.52) node   [align=left] {$\displaystyle \xi _{8}$};
\draw (506.33,495.67) node    {$\bullet $};
\draw (469,553.67) node    {$\bullet $};
\draw (544.33,554.67) node    {$\bullet $};

\end{tikzpicture}

\caption{$\mathbf{Q}(\xi) = (\xi_1\xi_2, \xi_2\xi_3, \xi_3\xi_1, \xi_2\xi_4,  
\xi_5\xi_6, 
\xi_6^2, \xi_7\xi_8, \xi_7\xi_9)$ with $d=9$ and $n=8$. Note that the left component is of Type (i) with a cycle of odd length, the middle component is of Type (ii), and the right component is of Type (iii). It is easy to check that $J_\mathbf{Q}(\xi;i)$ has the form (\ref{eq:linear_factor}) if and only if $i \in \{7,8,9\}$.}


  \label{fig:legal}
\end{figure}

As we have mentioned, the proof of Theorem~\ref{thm:Jacobian} is actually much stronger than the statement itself. For example, it indicates that the only nontrivial way for ``$J_{\mathbf{Q}}\equiv 0$'' to happen is that the corresponding directed graph contains a cycle of even length. What's more, the structure excavated can be used to prove the following result, which basically tells us that ``$J_\mathbf{Q}\not\equiv0$'' is downward hereditary in monomial cases. Such heredity will provide significant convenience in the proof of Theorem~\ref{main thm}, when we deal with the linear terms through certain inductive steps.

\begin{theorem}\label{thm:inherit}
    Let $\tilde{d}\geq n\geq2$, $\tilde{\mathbf{Q}}(\xi) = (\xi_{\lambda_{2m-1}} \xi_{\lambda_{2m}})_{m=1}^n$ be an $n$-tuple of quadratic monomials on $\mathbb{R}^{\tilde{d}}$, and $\mathbf{Q}'(\xi):= \tilde{\mathbf{Q}}\setminus\{ \xi_{\lambda_{2n-1}}\xi_{\lambda_{2n}} \} = (\xi_{\lambda_{2m-1}} \xi_{\lambda_{2m}})_{m=1}^{n-1}$ be an $n'$-tuple $(n' = n-1)$ of quadratic monomials on $\mathbb{R}^{d'}$ {\rm(}$d' \leq \tilde{d}${\rm)}. Suppose that $J_{\tilde{\mathbf{Q}}}(\xi;i_1,...,i_{\tilde{d}-n}) \not\equiv 0$ for some $i_1,...,i_{\tilde{d}-n}$, then we have $d' \geq n'$ and $J_{\mathbf{Q}'}(\xi;i_1',...,i_{d'-n'}')\not\equiv0$ for some $i_1',...,i_{d'-n'}'$.
\end{theorem}
\begin{proof}
    Let us first explain why $d' \geq n'$. Note that $n' = n-1$ and $d' = \tilde{d}$, $\tilde{d}-1$, or $\tilde{d}-2$. Since $\tilde{d} \geq n$, if $d' = \tilde{d}$ or $\tilde{d}-1$, then we have $d' \geq \tilde{d} - 1 \geq n-1 = n'$. Also, if $\tilde{d} \geq n+1$, then we still have $d' \geq \tilde{d} - 2 \geq n-1 = n'$, as desired. Thus, the case $d'<n'$ can only happen when $d' = \tilde{d} - 2$ and $\tilde{d}=n$, which means $\xi_{\lambda_{2n-1}}\not=\xi_{\lambda_{2n}}$ and they both do not appear in $\{ \xi_{\lambda_j} :j=1,...,2n-2  \}$. In this case, by expanding the determinant along the $\lambda_{2n-1}$$j$th and $\lambda_{2n}$$j$th columns, we immediately see that $J_{\tilde{\mathbf{Q}}}(\xi) \equiv 0$, which contradicts our assumption that $J_{\tilde{\mathbf{Q}}}(\xi) \not\equiv 0$. In other words, the case $d' < n'$ can never happen, and so we always have $d' \geq n'$, as desired. Now the expression $J_{\mathbf{Q}'}(\xi;i_1',...,i_{d'-n'}')$ makes sense, and it remains to show that it is nonzero for some $i_1',...,i_{d'-n'}'$.
    
    As in the proof of Theorem~\ref{thm:Jacobian}, the quadratic monomials in $\tilde{\mathbf{Q}}$ must be all distinct. And we can construct the associated graph $G = (V,E)$ and derive the structural properties in exactly the same way, except that $\mathbf{Q}$ and $d$ are replaced by $\tilde{\mathbf{Q}}$ and $\tilde{d}$ respectively. To avoid repetition, let us only briefly summarize what we get from the proof of Theorem~\ref{thm:Jacobian}. First, there is a one-to-one correspondence between nonzero terms (in the permutation expansion of $J_{\tilde{\mathbf{Q}}}(\xi;i_1,...,i_{\tilde{d}-n})$) and admissible pairs (init,\,ter), the existence of which is guaranteed by $J_{\tilde{\mathbf{Q}}}(\xi;i_1,...,i_{\tilde{d}-n}) \not\equiv 0$. Second, the components $G_k = (V_k, E_k)$ ($1\leq k \leq K$) can be divided into three types.    \begin{enumerate}[label={\rm (\roman*)}]
    \item No loop or tagged vertex ($\# V_k = \# E_k$), $G_k$ contains exactly one cycle of odd length.
    \item No tagged vertex, exactly one loop ($\# V_k = \# E_k$), $G_k$ is a tree upon deleting the loop.
    \item No loop, exactly one tagged vertex ($\# V_k = \# E_k + 1$), $G_k$ is a tree.
\end{enumerate}
Finally, for any graph (possibly) with loops and tagged vertices, as long as each component belongs to one of the three types above, the corresponding Jacobian must be nonzero.

With these structural properties in hand, we start the proof of Theorem~\ref{thm:inherit}. The key observation is that the graph of $\mathbf{Q}'$ is obtained by deleting the edge $\xi_{\lambda_{2n-1}}\xi_{\lambda_{2n}}$ from $G$ (deleting an edge also deletes any vertex that is incident only to that edge). Suppose $\xi_{\lambda_{2n-1}}\xi_{\lambda_{2n}} \in G_k$, then such deletion may destroy the nice structure of $G_k$ above, but will not affect all the other components. Therefore, to prove $J_{\mathbf{Q}'}(\xi;i_1',...,i_{d'-n'}')\not\equiv0$ for some $i_1',...,i_{d'-n'}'$, it suffices to verify that (each component of) $H_k := G_k - \xi_{\lambda_{2n-1}}\xi_{\lambda_{2n}}$ will still be of one of the three types, up to a possible adjustment of the tagged vertices. There are many cases, and we will discuss them one by one.

Let us first assume $d' = \tilde{d}$. In this case, $d' - n' = \tilde{d} - n + 1$, which means we need to tag an additional vertex.

If $G_k$ is of Type (i), then let $C$ be the unique cycle. If $\xi_{\lambda_{2n-1}}\xi_{\lambda_{2n}} \in C$, then $H_k$ is a tree, and will be of Type (iii) once we tag an arbitrary vertex of it. If $\xi_{\lambda_{2n-1}}\xi_{\lambda_{2n}} \notin C$, then $H_k$ splits into two components $G_k'$ and $G_k''$. Without loss of generality, assume $C \subset G_k'$, then clearly $G_k'$ is of Type (i). On the other hand, $G_k''$ is a tree, and will be of Type (iii) once we tag an arbitrary vertex of it. This completes the proof when $G_k$ is of Type (i).

If $G_k$ is of Type (ii), then let $\xi_j^2$ be the unique loop. If $\xi_{\lambda_{2n-1}}\xi_{\lambda_{2n}} = \xi_j^2$, then $H_k$ is a tree, and will be of Type (iii) once we tag an arbitrary vertex of it. If $\xi_{\lambda_{2n-1}}\xi_{\lambda_{2n}} \neq \xi_j^2$, then $H_k$ splits into two components $G_k'$ and $G_k''$. Without loss of generality, assume $\xi_j^2 \in G_k'$, then clearly $G_k'$ is of Type (ii). On the other hand, $G_k''$ is a tree, and will be of Type (iii) once we tag an arbitrary vertex of it. This completes the proof when $G_k$ is of Type (ii).

If $G_k$ is of Type (iii), then let $\xi_j$ be the unique tagged vertex. Note that $H_k$ must split into two components $G_k'$ and $G_k''$. Without loss of generality, assume $\xi_j \in G_k'$, then clearly $G_k'$ is of Type (iii). On the other hand, $G_k''$ is a tree, and will also be of Type (iii) once we tag an arbitrary vertex of it. This completes the proof when $G_k$ is of Type (iii).

Next, let us assume $d' = \tilde{d} - 1$. In this case, if $\xi_{\lambda_{2n-1}}=\xi_{\lambda_{2n}}$, then $G_k$ is a single loop, that is, $V_k = \{\xi_{\lambda_{2n}}\}$ and $E_k = \{\xi_{\lambda_{2n}}^2\}$, thus $H_k = \varnothing$ and the conclusion trivially holds true. Otherwise, $\xi_{\lambda_{2n-1}} \neq \xi_{\lambda_{2n}}$, which means $\xi_{\lambda_{2n-1}}\xi_{\lambda_{2n}}$ must be a pendant edge, and $H_k$ remains nonempty and connected. Also, $d' - n' = \tilde{d} - n$, which means we are only allowed to rearrange the tagged vertices without changing the number.

If $G_k$ is of Type (i), then since $\xi_{\lambda_{2n-1}}\xi_{\lambda_{2n}}$ is a pendant edge that does not affect the cycle, $H_k$ is still of Type (i). Similarly, if $G_k$ is of Type (ii), then since $\xi_{\lambda_{2n-1}}\xi_{\lambda_{2n}}$ is a pendant edge that is not the loop, $H_k$ is still of Type (ii). 

If $G_k$ is of Type (iii), then let $\xi_j$ be the unique tagged vertex. Without loss of generality, let $\xi_{\lambda_{2n}}$ be the (only) pendant vertex of $\xi_{\lambda_{2n-1}}\xi_{\lambda_{2n}}$. If $\xi_j = \xi_{\lambda_{2n}}$, then $H_k$ is a tree, and will be of Type (iii) once we tag an arbitrary vertex of it. If $\xi_j \neq \xi_{\lambda_{2n}}$, then $H_k$ is a tree containing $\xi_j$, and is automatically of Type (iii).

Finally, let us assume $d' = \tilde{d}-2$. In this case, $d' - n' = \tilde{d} - n - 1$, which means we need to delete a tagged vertex.(Recall that we have shown $d' \geq n'$ at the beginning, so $\tilde{d} \geq n + 1$.) Besides, both $\xi_{\lambda_{2n-1}}$ and $\xi_{\lambda_{2n}}$ do not appear in $\{\xi_{\lambda_j}: j=1,...,2n-2\}$. Since $\xi_{\lambda_{2n-1}}\xi_{\lambda_{2n}} \in G_k$, which is connected, we must have either $V_k = \{\xi_{\lambda_{2n-1}}, \xi_{\lambda_{2n}}\}$ when $\xi_{\lambda_{2n-1}} \neq \xi_{\lambda_{2n}}$, or $V_k = \{\xi_{\lambda_{2n}}\}$ when $\xi_{\lambda_{2n-1}}=\xi_{\lambda_{2n}}$. The latter case cannot happen, because it means $d' = \tilde{d} - 1$, while we need $d' = \tilde{d} - 2$. Therefore, it suffices to consider the former case. 

Clearly, $G_k$ cannot be of Type (i), which requires $\# V_k \geq 3$. If $G_k$ is of Type (ii), then $E_k = \{\xi_{\lambda_{2n-1}}\xi_{\lambda_{2n}}, \xi_{\lambda_{2n-1}}^2\}$ or $\{\xi_{\lambda_{2n-1}}\xi_{\lambda_{2n}}, \xi_{\lambda_{2n}}^2\}$. So, $H_k$ is a single loop at $\xi_{\lambda_{2n-1}}$ or $\xi_{\lambda_{2n}}$, which means $d' = \tilde{d} - 1$, contradicting $d' = \tilde{d} - 2$ again. Therefore, $G_k$ must be of Type (iii) and $E_k = \{\xi_{\lambda_{2n-1}}\xi_{\lambda_{2n}}\}$, thus $H_k = \varnothing$ and the conclusion trivially holds true. This completes the proof of Theorem~\ref{thm:inherit}.
\end{proof}

Before ending this section, let us see what Theorem~\ref{thm:inherit} does for the legal example in Figure~\ref{fig:legal}. Remember that we are deleting an edge from the graph. Up to isomorphism, there are six possibilities, which we now list as follows (Figures~\ref{fig:13}--\ref{fig:79}).

\begin{figure}[htbp]
  \centering
  \hspace*{0.5cm}  

\tikzset{every picture/.style={line width=0.75pt}} 

\begin{tikzpicture}[x=0.75pt,y=0.75pt,yscale=-1,xscale=1]

\draw    (346.33,462.49) -- (298.33,501.16) ;
\draw    (439.1,470.68) -- (402.33,498.82) ;
\draw   (418.86,450.92) .. controls (418.86,440) and (427.92,431.16) .. (439.1,431.16) .. controls (450.27,431.16) and (459.33,440) .. (459.33,450.92) .. controls (459.33,461.83) and (450.27,470.68) .. (439.1,470.68) .. controls (427.92,470.68) and (418.86,461.83) .. (418.86,450.92) -- cycle ;
\draw    (557.33,443.67) -- (520,501.67) ;
\draw    (595.33,502.67) -- (520,501.67) ;
\draw    (264,445.16) -- (299,501.16) ;
\draw    (229,501.16) -- (299,501.16) ;

\draw (276.92,436.39) node   [align=left] {$\displaystyle \xi _{3}$};
\draw (222.92,509.39) node   [align=left] {$\displaystyle \xi _{1}$};
\draw (309.38,508.34) node   [align=left] {$\displaystyle \xi _{2}$};
\draw (345.47,447.86) node   [align=left] {$\displaystyle \xi _{4}$};
\draw (264,445.16) node    {$\bullet $};
\draw (229,501.16) node    {$\bullet $};
\draw (299,501.16) node    {$\bullet $};
\draw (346.33,462.49) node    {$\bullet $};
\draw (411.71,508.01) node   [align=left] {$\displaystyle \xi _{5}$};
\draw (447.47,480.86) node   [align=left] {$\displaystyle \xi _{6}$};
\draw (439.1,470.68) node    {$\bullet $};
\draw (402.33,498.82) node    {$\bullet $};
\draw (569.92,434.57) node   [align=left] {$\displaystyle \xi _{9}$};
\draw (510.92,509.57) node   [align=left] {$\displaystyle \xi _{7}$};
\draw (606.38,510.52) node   [align=left] {$\displaystyle \xi _{8}$};
\draw (557.33,443.67) node    {$\bullet $};
\draw (520,501.67) node    {$\bullet $};
\draw (595.33,502.67) node    {$\bullet $};

\end{tikzpicture}

\caption{$\mathbf{Q}'(\xi) = (\xi_1\xi_2, \xi_2\xi_3, \xi_2\xi_4,  
\xi_5\xi_6, 
\xi_6^2, \xi_7\xi_8, \xi_7\xi_9)$ with $d'=9$ and $n'=7$. The left component becomes Type (iii), and $J_{\mathbf{Q}'}(\xi;i_1',i_2') \not\equiv 0$ if and only if $i_1' \in \{1,2,3,4\}$ and $i_2' \in \{7,8,9\}$, or $i_1' \in \{7,8,9\}$ and $i_2' \in \{1,2,3,4\}$.}


  \label{fig:13}
\end{figure}
\begin{figure}[htbp]
  \centering
  \hspace*{0.5cm}  

\tikzset{every picture/.style={line width=0.75pt}} 

\begin{tikzpicture}[x=0.75pt,y=0.75pt,yscale=-1,xscale=1]

\draw    (352.33,610.49) -- (304.33,649.16) ;
\draw    (445.1,618.68) -- (408.33,646.82) ;
\draw   (424.86,598.92) .. controls (424.86,588) and (433.92,579.16) .. (445.1,579.16) .. controls (456.27,579.16) and (465.33,588) .. (465.33,598.92) .. controls (465.33,609.83) and (456.27,618.68) .. (445.1,618.68) .. controls (433.92,618.68) and (424.86,609.83) .. (424.86,598.92) -- cycle ;
\draw    (563.33,591.67) -- (526,649.67) ;
\draw    (601.33,650.67) -- (526,649.67) ;
\draw    (270,593.16) -- (304.33,649.16) ;
\draw    (235,649.16) -- (270,593.16) ;

\draw (282.92,584.39) node   [align=left] {$\displaystyle \xi _{3}$};
\draw (228.92,657.39) node   [align=left] {$\displaystyle \xi _{1}$};
\draw (315.38,656.34) node   [align=left] {$\displaystyle \xi _{2}$};
\draw (351.47,595.86) node   [align=left] {$\displaystyle \xi _{4}$};
\draw (270,593.16) node    {$\bullet $};
\draw (235,649.16) node    {$\bullet $};
\draw (304.33,649.16) node    {$\bullet $};
\draw (352.33,610.49) node    {$\bullet $};
\draw (417.71,656.01) node   [align=left] {$\displaystyle \xi _{5}$};
\draw (453.47,628.86) node   [align=left] {$\displaystyle \xi _{6}$};
\draw (445.1,618.68) node    {$\bullet $};
\draw (408.33,646.82) node    {$\bullet $};
\draw (575.92,582.57) node   [align=left] {$\displaystyle \xi _{9}$};
\draw (516.92,657.57) node   [align=left] {$\displaystyle \xi _{7}$};
\draw (612.38,658.52) node   [align=left] {$\displaystyle \xi _{8}$};
\draw (563.33,591.67) node    {$\bullet $};
\draw (526,649.67) node    {$\bullet $};
\draw (601.33,650.67) node    {$\bullet $};

\end{tikzpicture}

\caption{$\mathbf{Q}'(\xi) = (\xi_1\xi_3, \xi_2\xi_3, \xi_2\xi_4,  
\xi_5\xi_6, 
\xi_6^2, \xi_7\xi_8, \xi_7\xi_9)$ with $d'=9$ and $n'=7$. The analysis is the same as that for Figure~\ref{fig:13}.}


  \label{fig:12}
\end{figure}
\begin{figure}[htbp]
  \centering
  \hspace*{0.5cm}  

\tikzset{every picture/.style={line width=0.75pt}} 

\begin{tikzpicture}[x=0.75pt,y=0.75pt,yscale=-1,xscale=1]

\draw   (267,739.55) -- (302,795.55) -- (232,795.55) -- cycle ;
\draw    (442.1,765.07) -- (405.33,793.22) ;
\draw   (421.86,745.31) .. controls (421.86,734.4) and (430.92,725.55) .. (442.1,725.55) .. controls (453.27,725.55) and (462.33,734.4) .. (462.33,745.31) .. controls (462.33,756.23) and (453.27,765.07) .. (442.1,765.07) .. controls (430.92,765.07) and (421.86,756.23) .. (421.86,745.31) -- cycle ;
\draw    (560.33,738.06) -- (523,796.06) ;
\draw    (598.33,797.06) -- (523,796.06) ;

\draw (279.92,730.79) node   [align=left] {$\displaystyle \xi _{3}$};
\draw (225.92,803.79) node   [align=left] {$\displaystyle \xi _{1}$};
\draw (312.38,802.74) node   [align=left] {$\displaystyle \xi _{2}$};
\draw (267,739.55) node    {$\bullet $};
\draw (232,795.55) node    {$\bullet $};
\draw (302,795.55) node    {$\bullet $};
\draw (414.71,802.41) node   [align=left] {$\displaystyle \xi _{5}$};
\draw (450.47,775.26) node   [align=left] {$\displaystyle \xi _{6}$};
\draw (442.1,765.07) node    {$\bullet $};
\draw (405.33,793.22) node    {$\bullet $};
\draw (572.92,728.97) node   [align=left] {$\displaystyle \xi _{9}$};
\draw (513.92,803.97) node   [align=left] {$\displaystyle \xi _{7}$};
\draw (609.38,804.92) node   [align=left] {$\displaystyle \xi _{8}$};
\draw (560.33,738.06) node    {$\bullet $};
\draw (523,796.06) node    {$\bullet $};
\draw (598.33,797.06) node    {$\bullet $};

\end{tikzpicture}

\caption{$\mathbf{Q}'(\xi) = (
\xi_1\xi_2, \xi_2\xi_3, 
\xi_3\xi_1,  
\xi_5\xi_6, 
\xi_6^2, \xi_7\xi_8, \xi_7\xi_9)$ with $d'=8$ and $n'=7$. The left component remains Type (i), and $J_{\mathbf{Q}'}(\xi;i')\not\equiv 0$ if and only if $i' \in \{7,8,9\}$.}


  \label{fig:24}
\end{figure}
\begin{figure}[htbp]
  \centering
  \hspace*{0.5cm}  

\tikzset{every picture/.style={line width=0.75pt}} 

\begin{tikzpicture}[x=0.75pt,y=0.75pt,yscale=-1,xscale=1]

\draw   (269,884.16) -- (304,940.16) -- (234,940.16) -- cycle ;
\draw    (351.33,901.49) -- (303.33,940.16) ;
\draw   (423.86,889.92) .. controls (423.86,879) and (432.92,870.16) .. (444.1,870.16) .. controls (455.27,870.16) and (464.33,879) .. (464.33,889.92) .. controls (464.33,900.83) and (455.27,909.68) .. (444.1,909.68) .. controls (432.92,909.68) and (423.86,900.83) .. (423.86,889.92) -- cycle ;
\draw    (562.33,882.67) -- (525,940.67) ;
\draw    (600.33,941.67) -- (525,940.67) ;

\draw (281.92,875.39) node   [align=left] {$\displaystyle \xi _{3}$};
\draw (227.92,948.39) node   [align=left] {$\displaystyle \xi _{1}$};
\draw (314.38,947.34) node   [align=left] {$\displaystyle \xi _{2}$};
\draw (350.47,886.86) node   [align=left] {$\displaystyle \xi _{4}$};
\draw (269,884.16) node    {$\bullet $};
\draw (234,940.16) node    {$\bullet $};
\draw (304,940.16) node    {$\bullet $};
\draw (351.33,901.49) node    {$\bullet $};
\draw (452.47,919.86) node   [align=left] {$\displaystyle \xi _{6}$};
\draw (444.1,909.68) node    {$\bullet $};
\draw (574.92,873.57) node   [align=left] {$\displaystyle \xi _{9}$};
\draw (515.92,948.57) node   [align=left] {$\displaystyle \xi _{7}$};
\draw (611.38,949.52) node   [align=left] {$\displaystyle \xi _{8}$};
\draw (562.33,882.67) node    {$\bullet $};
\draw (525,940.67) node    {$\bullet $};
\draw (600.33,941.67) node    {$\bullet $};

\end{tikzpicture}

\caption{$\mathbf{Q}'(\xi) = (\xi_1\xi_2, \xi_2\xi_3,
\xi_3\xi_1, \xi_2\xi_4,  
\xi_6^2, \xi_7\xi_8, \xi_7\xi_9)$ with $d'=8$ and $n'=7$. The middle component remains Type (ii), and $J_{\mathbf{Q}'}(\xi;i')\not\equiv 0$ if and only if $i' \in \{7,8,9\}$.}


  \label{fig:56}
\end{figure}
\begin{figure}[htbp]
  \centering
  \hspace*{0.5cm}  

\tikzset{every picture/.style={line width=0.75pt}} 

\begin{tikzpicture}[x=0.75pt,y=0.75pt,yscale=-1,xscale=1]

\draw   (269,1029.22) -- (304,1085.22) -- (234,1085.22) -- cycle ;
\draw    (351.33,1046.55) -- (303.33,1085.22) ;
\draw    (444.1,1054.74) -- (407.33,1082.89) ;
\draw    (562.33,1027.73) -- (525,1085.73) ;
\draw    (600.33,1086.73) -- (525,1085.73) ;

\draw (281.92,1020.46) node   [align=left] {$\displaystyle \xi _{3}$};
\draw (227.92,1093.46) node   [align=left] {$\displaystyle \xi _{1}$};
\draw (314.38,1092.41) node   [align=left] {$\displaystyle \xi _{2}$};
\draw (350.47,1031.93) node   [align=left] {$\displaystyle \xi _{4}$};
\draw (269,1029.22) node    {$\bullet $};
\draw (234,1085.22) node    {$\bullet $};
\draw (304,1085.22) node    {$\bullet $};
\draw (351.33,1046.55) node    {$\bullet $};
\draw (416.71,1092.07) node   [align=left] {$\displaystyle \xi _{5}$};
\draw (452.47,1064.93) node   [align=left] {$\displaystyle \xi _{6}$};
\draw (444.1,1054.74) node    {$\bullet $};
\draw (407.33,1082.89) node    {$\bullet $};
\draw (574.92,1018.64) node   [align=left] {$\displaystyle \xi _{9}$};
\draw (515.92,1093.64) node   [align=left] {$\displaystyle \xi _{7}$};
\draw (611.38,1094.58) node   [align=left] {$\displaystyle \xi _{8}$};
\draw (562.33,1027.73) node    {$\bullet $};
\draw (525,1085.73) node    {$\bullet $};
\draw (600.33,1086.73) node    {$\bullet $};

\end{tikzpicture}

\caption{$\mathbf{Q}'(\xi) = (\xi_1\xi_2, \xi_2\xi_3,
\xi_3\xi_1, \xi_2\xi_4,  
\xi_5\xi_6, \xi_7\xi_8, \xi_7\xi_9)$ with $d'=9$ and $n'=7$. The middle component becomes Type (iii), and $J_{\mathbf{Q}'}(\xi;i_1',i_2')\not\equiv 0$ if and only if $i_1' \in \{5,6\}$ and $i_2' \in \{7,8,9\}$, or $i_1' \in \{7,8,9\}$ and $i_2' \in \{5,6\}$.}


  \label{fig:66}
\end{figure}
\begin{figure}[htbp]
  \centering
  \hspace*{0.5cm}  

\tikzset{every picture/.style={line width=0.75pt}} 

\begin{tikzpicture}[x=0.75pt,y=0.75pt,yscale=-1,xscale=1]

\draw   (269,1178.22) -- (304,1234.22) -- (234,1234.22) -- cycle ;
\draw    (351.33,1195.55) -- (303.33,1234.22) ;
\draw    (444.1,1203.74) -- (407.33,1231.89) ;
\draw   (423.86,1183.98) .. controls (423.86,1173.07) and (432.92,1164.22) .. (444.1,1164.22) .. controls (455.27,1164.22) and (464.33,1173.07) .. (464.33,1183.98) .. controls (464.33,1194.89) and (455.27,1203.74) .. (444.1,1203.74) .. controls (432.92,1203.74) and (423.86,1194.89) .. (423.86,1183.98) -- cycle ;
\draw    (600.33,1235.73) -- (525,1234.73) ;

\draw (281.92,1169.46) node   [align=left] {$\displaystyle \xi _{3}$};
\draw (227.92,1242.46) node   [align=left] {$\displaystyle \xi _{1}$};
\draw (314.38,1241.41) node   [align=left] {$\displaystyle \xi _{2}$};
\draw (350.47,1180.93) node   [align=left] {$\displaystyle \xi _{4}$};
\draw (269,1178.22) node    {$\bullet $};
\draw (234,1234.22) node    {$\bullet $};
\draw (304,1234.22) node    {$\bullet $};
\draw (351.33,1195.55) node    {$\bullet $};
\draw (416.71,1241.07) node   [align=left] {$\displaystyle \xi _{5}$};
\draw (452.47,1213.93) node   [align=left] {$\displaystyle \xi _{6}$};
\draw (444.1,1203.74) node    {$\bullet $};
\draw (407.33,1231.89) node    {$\bullet $};
\draw (515.92,1242.64) node   [align=left] {$\displaystyle \xi _{7}$};
\draw (611.38,1243.58) node   [align=left] {$\displaystyle \xi _{8}$};
\draw (525,1234.73) node    {$\bullet $};
\draw (600.33,1235.73) node    {$\bullet $};

\end{tikzpicture}

\caption{$\mathbf{Q}'(\xi) = (\xi_1\xi_2, \xi_2\xi_3,
\xi_3\xi_1, \xi_2\xi_4,  
\xi_5\xi_6,
\xi_6^2, \xi_7\xi_8)$ with $d'=8$ and $n'=7$. The right component remains Type (iii), and $J_{\mathbf{Q}'}(\xi;i')\not\equiv 0$ if and only if $i' \in \{7,8\}$.}


  \label{fig:79}
\end{figure}

\section{Broad--narrow analysis and iteration}\label{section3}

In this section, we introduce the broad--narrow analysis and its iteration adapted to a general quadratic form. Besides, to prepare for the proof of Theorem~\ref{main thm} in the next section, we also establish the  main tools that will be needed to deal with the broad part and the narrow part: For the broad part, we prove bilinear restriction estimates independent of the transversality parameters; while for the narrow part, we recall classical decoupling results.

If $\mathbf{Q}$ satisfies the constraints of Theorem \ref{main thm}, then $J_{\mathbf{Q}}(\xi;i_1,...,i_{d-n}) \not\equiv 0$ for some $i_1,$...,$i_{d-n}$. By Theorem~\ref{thm:Jacobian}, there exist $i_1,...,i_{d-n}$ (possibly different) such that
\begin{equation*}
   	|J_{\mathbf{Q}}(\xi;i_1,...,i_{d-n})|\sim \prod_{j=1}^d |\xi_{j}|^{s_{j}}
\end{equation*}
with $s_j\in\N$ {\rm(}$1\leq j\leq d${\rm)}, $\sum_{j=1}^d s_{j}=n$, and $\max_j s_j \leq \max_j w_j - 1$. After ignoring those $s_j$ with zero value and then reindexing the variables (this does not affect the boundedness of restriction estimates), we may without loss of generality assume 
\begin{equation}\label{s2s3e1}
    |J_{\mathbf{Q}}(\xi;i_1,...,i_{d-n})|\sim \prod_{j=1}^t |\xi_{j}|^{s_{j}}
\end{equation}
with $s_j\in\N^+$ {\rm(}$1\leq j\leq t${\rm)}, $\sum_{j=1}^t s_{j}=n$, and $\max_j s_j \leq \max_j w_j - 1$. This linear factorization is crucial in our argument, since it will naturally induce an appropriate transversality condition and thus the associated broad--narrow analysis.

\subsection{Broad--narrow analysis and iteration}\label{subsection31}
From now on, we consider a large scale $R\geq 1$. Let us start by defining several constants.

\begin{definition}\label{s2def1}
    Let $R\geq1$ and $1\leq \mu_1,...,\mu_d \leq R^{1/2}$. Suppose that $\mathbf{Q}=(Q_1,...,Q_n)$ is an $n$-tuple of quadratic forms in $d$ variables. 
    
    Define $D_p(\mu_1,...,\mu_d;R)$ to be the smallest constant such that
    \begin{equation*}
        \|E^{\mathbf{Q}}f\|_{L^p(B^{d+n}_R)} \leq D_p(\mu_1,\dots,\mu_d;R) \|f\|_{L^p(\mathbb{R}^d)}
    \end{equation*}
    holds for all $f$ satisfying $\on{supp} f\subset [a_1,a_1+\mu_1^{-1}]\times\cdots \times [a_d,a_d+\mu_d^{-1}]\subset [0,1]^d$ for some $a_1,...,a_d$. When $\mu_1=\cdots=\mu_d=1$, we abbreviate $D_p(1,...,1;R)$ as $D_p(R)$. 
    
    Define $BD_p(\mu_1,...,\mu_d;R)$ to be the smallest constant such that
    \begin{equation*}
        \Big\||E^{\mathbf{Q}}f_1 E^{\mathbf{Q}}f_2|^{\frac{1}{2}}\Big\|_{L^p(B^{d+n}_R)} \leq BD_p(\mu_1,...,\mu_d;R) \|f_1\|_{L^p(\mathbb{R}^d)}^{\frac{1}{2}} \|f_2\|_{L^p(\mathbb{R}^d)}^{\frac{1}{2}}
    \end{equation*}
    holds for all $f_1$ satisfying $\on{supp} f_1 \subset [a_1,a_1+\mu_1^{-1}] \times\cdots\times [a_d,a_d+\mu_d^{-1}]\subset [0,1]^d$ and  $f_2$ satisfying $\on{supp} f_2 \subset [b_1,b_1+\mu_1^{-1}] \times\cdots\times [b_d,b_d+\mu_d^{-1}]\subset [0,1]^d$ with $|a_j-b_j|\geq 10\mu_j^{-1}$ whenever $\mu_j>10$ \rm{(}$1\leq j \leq d$\rm{)}. If $1<\mu_j\leq 10$ \rm{(}$1\leq j \leq d$\rm{)}, we say $f_1$ and $f_2$ are both supported in $[0,1]$ with respect to the $j$th coordinate.  
\end{definition}
Note that 
\begin{equation}\label{s2e1}
    D_p(\mu_1,...,\mu_d;R) \leq  D_p(\mu'_1,...,\mu'_d;R)
\end{equation}
whenever $\mu_j \geq \mu_j'$ for all $1\leq j\leq d$.

For the convenience of discussion, we will only provide details when $t=2$ in (\ref{s2s3e1}), that is,
\begin{equation*}
    |J_{\mathbf{Q}}(\xi;i_1,...,i_{d-n})|\sim |\xi_{1}|^{s_{1}}|\xi_{2}|^{s_{2}},
\end{equation*}
with $s_1,s_2 \in \N^+$ and $s_1+s_2=n$. General cases can be handled similarly.

Now we start to perform the bilinear broad--narrow analysis for $\mathbf{Q}$ with $t=2$. 
Let $K$ be a large dyadic integer satisfying $K \sim \log R$. Divide $[0,1]^d$ into cubes $\tau$ of side length $K^{-1}$, and  write $f=\sum_\tau f_\tau$. For each $x\in B_R$, we define its  \textit{significant set} as
$$    \mathcal{S}(x) :=\Big\{ \tau:  |E^{\mathbf{Q}} f _\tau(x) |\geq \frac{1}{100\# \{\tau\}}   |E^{\mathbf{Q}} f(x) |\Big\}.  $$
We say $x$ is \textit{broad} if there exist $\tau_1,\tau_2 \in \mathcal{S}(x)$ such that for any $\xi=(\xi_1,...,\xi_d) \in \tau_1$, $\eta=(\eta_1,...,\eta_d) \in \tau_2$, we have
\begin{equation}\label{s2e6}
    |\xi_1-\eta_1|\geq 10K^{-1}, ~ |\xi_2-\eta_2|\geq 10K^{-1} .
\end{equation}
Otherwise, we say $x$ is \textit{narrow}. If $x$ is broad, then by the definition of $\mathcal{S}(x)$, there exist $\tau_1$, $\tau_2$ satisfying (\ref{s2e6}) such that
$$  \lvert E^{\mathbf{Q}}f(x)\rvert \lesssim K^{d}   \Big| E^{\mathbf{Q}}f_{\tau_1}(x) E^{\mathbf{Q}}f_{\tau_2}(x)\Big|^{\frac{1}{2}} . $$
If $x$ is narrow, then there is
$$    \lvert E^{\mathbf{Q}}f(x)\rvert \lesssim  \Big|\sum_{\tau \subset W} E^{\mathbf{Q}}f_\tau(x)\Big| , $$
where $W$ denotes a rectangular box of size $\sim K^{-1}\times 1\times \cdots \times 1$ or $1 \times K^{-1}\times 1\times \cdots\times 1$. In summary, for general $x\in B_R$, there exists a universal constant $C$ (such as $100$) so that
\begin{align}
    \lvert E^{\mathbf{Q}}f(x)\rvert \leq &CK^{d} \max_{\tau_1,\tau_2 \text{~satisfy~} (\ref{s2e6})} | E^{\mathbf{Q}}f_{\tau_1}(x) E^{\mathbf{Q}}f_{\tau_2}(x)|^{\frac{1}{2}}  \nonumber\\
    &+C\max_{W_{1,0}}\Big| E^{\mathbf{Q}}f_{W_{1,0}}(x)\Big| +C\max_{W_{0,1}}\Big| E^{\mathbf{Q}}f_{W_{0,1}}(x)\Big|,   \label{s2s3e3}
\end{align}
where $W_{1,0}$ and $W_{0,1}$ denote rectangular boxes of sizes $ K^{-1}\times 1\times \cdots \times 1$ and $1 \times K^{-1}\times 1\times \cdots\times 1$, respectively. We regard (\ref{s2s3e3}) as the first step of the iteration. 

To continue the iteration, let us first consider $E^{\mathbf{Q}}f_{W_{1,0}}$. Divide  $W_{1,0}$ into rectangular boxes $\tilde{\tau}$ of size $K^{-2}\times K^{-1}\times \cdots\times K^{-1}$ along its axes, and write $f_{W_{1,0}}=\sum_{\tilde{\tau}}f_{\tilde{\tau}}$. Note that 
$$   \#\{\tilde{\tau}:\tilde{\tau}\subset W_{1,0}   \}\sim K^d.   $$
Similarly, for each $x\in B_R$, we define the significant set in $W_{1,0}$ as
$$    \tilde{\mathcal{S}}(x) :=\Big\{ \tilde{\tau}:  |E^{\mathbf{Q}} f_{\tilde{\tau}}(x) |\geq \frac{1}{100\# \{\tilde{\tau}\}}   |E^{\mathbf{Q}} f_{W_{1,0}}(x) |\Big\}.  $$
We say $x$ is broad if there exist $\tilde{\tau}_1,\tilde{\tau}_2 \in \tilde{\mathcal{S}}(x)$ such that for any $\xi=(\xi_1,...,\xi_d) \in \tilde{\tau}_1$, $\eta=(\eta_1,...,\eta_d) \in \tilde{\tau}_2$, we have
\begin{equation}\label{s2e60}
    |\xi_1-\eta_1|\geq 10K^{-2}, ~ |\xi_2-\eta_2|\geq 10K^{-1} .
\end{equation}
Otherwise, we say $x$ is narrow. Via a similar broad--narrow analysis as in the first step, we get
\begin{align}
    \lvert E^{\mathbf{Q}}f_{W_{1,0}}(x)\rvert \leq &CK^{d} \max_{\tilde{\tau}_1,\tilde{\tau}_2 \text{~satisfy~} (\ref{s2e60})} | E^{\mathbf{Q}}f_{\tilde{\tau}_1}(x) E^{\mathbf{Q}}f_{\tilde{\tau}_2}(x)|^{\frac{1}{2}}  \nonumber\\
    &+C\max_{W_{2,0}}\Big| E^{\mathbf{Q}}f_{W_{2,0}}(x)\Big| +C\max_{W_{1,1}}\Big| E^{\mathbf{Q}}f_{W_{1,1}}(x)\Big|.   \label{s2s3e4}
\end{align}
where $W_{2,0}$ and $W_{1,1}$ denote rectangular boxes of sizes $  K^{-2}\times1\times \cdots \times 1$ and $K^{-1} \times K^{-1}\times 1\times \cdots\times 1$, respectively.
Similarly, for $E^{\mathbf{Q}}f_{W_{0,1}}$, we can bound it by one broad term (bilinear term) and two narrow terms ($ E^{\mathbf{Q}}f_{W_{0,2}}$ and $ E^{\mathbf{Q}}f_{W_{1,1}}$). Plugging these estimates into (\ref{s2s3e3}), we get
\begin{align}
    \lvert E^{\mathbf{Q}}f(x)\rvert \leq &   3C^2 K^{d} \max_{(\ast)} | E^{\mathbf{Q}}f_{\tau_1}(x) E^{\mathbf{Q}}f_{\tau_2}(x)|^{\frac{1}{2}}  \label{s2s3e5}\\
    &+C^2\max_{W_{2,0}}\Big| E^{\mathbf{Q}}f_{W_{2,0}}(x)\Big| +2C^2\max_{W_{1,1}}\Big| E^{\mathbf{Q}}f_{W_{1,1}}(x)\Big| +C^2\max_{W_{0,2}}\Big| E^{\mathbf{Q}}f_{W_{0,2}}(x)\Big|.   \label{s2s3e6}
\end{align}
Let us explain the meaning of $(\ast)$ in \eqref{s2s3e5}. There are three possibilities for $\tau_1$ and $\tau_2$ in the broad term:

\noindent (1) $\tau_1$ and $\tau_2$ both have size $K^{-1}\times \cdots \times K^{-1}$, and for any $\xi=(\xi_1,...,\xi_d) \in \tau_1$, $\eta=(\eta_1,...,\eta_d) \in \tau_2$, we have
$$       |\xi_1-\eta_1|\geq 10K^{-1}, ~ |\xi_2-\eta_2|\geq 10K^{-1} ;   $$

\noindent (2) $\tau_1$ and $\tau_2$ both have size $K^{-2}\times K^{-1}\times\cdots \times K^{-1}$, and for any $\xi=(\xi_1,...,\xi_d) \in \tau_1$, $\eta=(\eta_1,...,\eta_d) \in \tau_2$, we have
$$       |\xi_1-\eta_1|\geq 10K^{-2}, ~ |\xi_2-\eta_2|\geq 10K^{-1} ;   $$

\noindent (3) $\tau_1$ and $\tau_2$ both have size $K^{-1}\times K^{-2}\times K^{-1}\times\cdots \times K^{-1}$, and for any $\xi=(\xi_1,...,\xi_d) \in \tau_1$, $\eta=(\eta_1,...,\eta_d) \in \tau_2$, we have
$$       |\xi_1-\eta_1|\geq 10K^{-1}, ~ |\xi_2-\eta_2|\geq 10K^{-2} .   $$
Besides, in either case, we always have $\#\{\tau_1\},\#\{\tau_2\}\sim K^d$. 
For the convenience of the following iterative discussions, we uniformly write them as: $\tau_1$ and $\tau_2$ denote rectangular boxes of the same size $\mu_1^{-1}\times \mu_2^{-1}\times 1\times \cdots \times 1$ (where the last $d-2$ components are originally $K^{-1}$ but can be harmlessly interpreted as $1$ by zero extensions of $f_{\tau_1}$ and $f_{\tau_2}$) for some $1<\mu_1\leq R^{1/2}$ and $1<\mu_2\leq R^{1/2}$, $\#\{\tau_1\},\#\{\tau_2\}\sim K^d$, and for any $\xi=(\xi_1,...,\xi_d) \in \tau_1$, $\eta=(\eta_1,...,\eta_d) \in \tau_2$, we have
\begin{equation}\label{s2s3e7}
    |\xi_1-\eta_1|\geq 10\mu_1^{-1}, ~ |\xi_2-\eta_2|\geq 10\mu_2^{-1} .
\end{equation}
This just is the meaning of $(\ast)$. 
Next, we apply the same type of broad--narrow analysis to the term $E^{\mathbf{Q}}f_{W_{1,1}}$ in (\ref{s2s3e6}) again. Then it is bounded by one broad term (bilinear term) and two narrow terms ($ E^{\mathbf{Q}}f_{W_{2,1}}$ and $ E^{\mathbf{Q}}f_{W_{1,2}}$). Thus, (\ref{s2s3e6}) further becomes 
\begin{align}
    \lvert E^{\mathbf{Q}}f(x)\rvert \leq & 5C^3 K^{d} \max_{(\ast)} | E^{\mathbf{Q}}f_{\tau_1}(x) E^{\mathbf{Q}}f_{\tau_2}(x)|^{\frac{1}{2}}  \nonumber\\
    &+2C^2\sum_{j=0}^1 \max_{W_{2,j}}\Big| E^{\mathbf{Q}}f_{W_{2,j}}(x)\Big| +2C^2\sum_{j=0}^1\max_{W_{j,2}}\Big| E^{\mathbf{Q}}f_{W_{j,2}}(x)\Big|,   \label{s2s3e8}
\end{align}
where $(\ast)$ has the same meaning as in (\ref{s2s3e5}), and $W_{2,j}$ and $W_{j,2}$ denote rectangular boxes of sizes $ K^{-2}\times K^{-j}\times1\times \cdots \times 1$ and $K^{-j} \times K^{-2}\times 1\times \cdots\times 1$, respectively. We regard (\ref{s2s3e8}) as the second step of the iteration.

Suppose that $K^m=R^{1/2}$ for some $m$. Repeating the above arguments until the $m$$j$th step, we obtain
\begin{align}
    \lvert E^{\mathbf{Q}}f(x)\rvert \lesssim & m^2 (4C^2)^{m} K^{d} \max_{(\ast)} | E^{\mathbf{Q}}f_{\tau_1}(x) E^{\mathbf{Q}}f_{\tau_2}(x)|^{\frac{1}{2}}  \nonumber\\
    &+(4C^2)^{m}\sum_{j=0}^{m-1}\max_{W_{m,j}}\Big| E^{\mathbf{Q}}f_{W_{m,j}}(x)\Big| +(4C^2)^{m}\sum_{j=0}^{m-1}\max_{W_{j,m}}\Big| E^{\mathbf{Q}}f_{W_{j,m}}(x)\Big|,   \label{s2s3e9}
\end{align}
where $(\ast)$ has the same meaning as in (\ref{s2s3e5}), and $W_{m,j}$ and $W_{j,m}$ denote rectangular boxes of sizes $ R^{-1/2}\times K^{-j}\times1\times \cdots \times 1$ and $K^{-j} \times R^{-1/2}\times 1\times \cdots\times 1$, respectively. 

We integrate over $B_R$ on both sides of (\ref{s2s3e9}), and use (\ref{s2e1}) to obtain
\begin{align*}
    \| E^{\mathbf{Q}}f\|_{L^p(B_R)} &\lesssim  m^2 (4C^2)^{m} K^{d} \sum_{(\ast)} \Big\|| E^{\mathbf{Q}}f_{\tau_1} E^{\mathbf{Q}}f_{\tau_2}|^{\frac{1}{2}}\Big\|_{L^p(B_R)}  \\
    &\quad\quad+ (4C^2)^{m}\sum_{j=0}^{m-1}\Big(\sum_{W_{m,j}}\| E^{\mathbf{Q}}f_{W_{m,j}}\|^p_{L^p(B_R)}\Big)^{\frac{1}{p}} + (4C^2)^{m}\sum_{j=0}^{m-1}\Big(\sum_{W_{j,m}}\| E^{\mathbf{Q}}f_{W_{j,m}}\|^p_{L^p(B_R)} \Big)^{\frac{1}{p}} \\ 
    &\lesssim m^2 (4C^2)^{m} K^{3d}\sup_{1<\mu_1,\mu_2\leq R^{1/2}} BD_p(\mu_1,\mu_2,1,...,1;R)\|f\|_{L^p}    \\
    &  \quad\quad + (4C^2)^{m} \sum_{j=0}^{m-1}\Big[\sum_{W_{m,j}}  (D_p(R^{\frac{1}{2}},K^j,1,...,1;R)\|f_{W_{m,j}}\|_{L^p})^p   \Big]^{\frac{1}{p}} \\
    &\quad\quad+ (4C^2)^{m}\sum_{j=0}^{m-1}\Big[\sum_{W_{j,m}}(D_p(K^j,R^{\frac{1}{2}},1,...,1;R) \|f_{W_{j,m}}\|_{L^p})^p    \Big]^{\frac{1}{p}}   \\
    &  \leq   m^2 (4C^2)^{m} K^{3d}\sup_{1<\mu_1,\mu_2\leq R^{1/2}} BD_p(\mu_1,\mu_2,1,...,1;R)\|f\|_{L^p}    \\
    &  \quad\quad + (4C^2)^{m} D_p(R^{\frac{1}{2}},1,...,1;R)\sum_{j=0}^{m-1}\Big(\sum_{W_{m,j}}  \|f_{W_{m,j}}\|_{L^p}^p   \Big)^{\frac{1}{p}} \\
    &\quad\quad+(4C^2)^{m}D_p(1,R^{\frac{1}{2}},1,...,1;R)\sum_{j=0}^{m-1}\Big(\sum_{W_{j,m}} \|f_{W_{j,m}}\|_{L^p}^p    \Big)^{\frac{1}{p}}   \\
    & = m^2 (4C^2)^{m} K^{3d}\sup_{1<\mu_1,\mu_2\leq R^{1/2}} BD_p(\mu_1,\mu_2,1,...,1;R)\|f\|_{L^p}    \\
    & \quad\quad + m (4C^2)^{m} D_p(R^{\frac{1}{2}},1,...,1;R)    \|f\|_{L^p}+ m (4C^2)^{m}   D_p(1,R^{\frac{1}{2}},1,...,1;R) \|f\|_{L^p}  \\
    & \lesssim R^{\epsilon} \sup_{1<\mu_1,\mu_2\leq R^{1/2}} BD_p(\mu_1,\mu_2,1,...,1;R)\|f\|_{L^p}    \\
    & \quad\quad + R^{\epsilon} D_p(R^{\frac{1}{2}},1,...,1;R)    \|f\|_{L^p}+ R^{\epsilon}   D_p(1,R^{\frac{1}{2}},1,...,1;R) \|f\|_{L^p},
\end{align*}
where $0<\epsilon\ll1 $ is any fixed small constant, and we used  $m^2 (4C^2)^{m} K^{3d} \lesssim R^{\epsilon}$ by noting $K\sim \log R$. 
Therefore, by the definition of $D_p(R)$, we finally get
\begin{equation*}
D_p(R) \lesssim R^{\epsilon}  \sup_{1<\mu_1,\mu_2\leq R^{1/2}} BD_p(\mu_1,\mu_2,1,...,1;R)+  R^{\epsilon} D_p(R^{\frac{1}{2}},1,...,1;R)+ R^{\epsilon} D_p(1,R^{\frac{1}{2}},1,...,1;R).
\end{equation*}

If we are in the general case (\ref{s2s3e1}) where $t$ could be larger than $2$, that is,
\begin{equation*}
    |J_{\mathbf{Q}}(\xi;i_1,...,i_{d-n})|\sim |\xi_{1}|^{s_{1}}|\xi_{2}|^{s_{2}}\cdots |\xi_t|^{s_t},
\end{equation*}
with $s_1,s_2,...,s_t \in \N^+$ and $s_1+s_2+\cdots +s_t=n$, then the corresponding iterative formula becomes
\begin{equation}\label{s2s3e10}
D_p(R) \lesssim R^{\epsilon} \sup_{1<\mu_1,...,\mu_t\leq R^{1/2}} BD_p(\mu_1,...,\mu_t,1,...,1;R)+  R^{\epsilon} \sum_{j=1}^t D_p(1,...,1,\underset{\substack{\uparrow \\ j}}{R^\frac{1}{2}},1,...,1;R).
\end{equation}
We omit the details of the proof of (\ref{s2s3e10}) as it is very similar to the $t=2$ case: Keep performing previous broad--narrow analysis to reduce original $f$ with frequency support $[0,1]^d$ to some $f_j$ whose frequency supports have dimensions   
 $$1\times \cdots\times 1\times \underset{\substack{\uparrow \\ j}}{R^{-\frac{1}{2}}}\times 1\times \cdots \times 1,\quad \quad \quad \quad 1\leq j \leq t.$$ 
Though the constants depend on $t$,  since $t\leq d$, they would still be eventually absorbed by $R^\epsilon$.
So far, we have completed the broad--narrow analysis, which reduces estimating $D_p(R)$ to handle the broad part and narrow part separately. 

\subsection{A refined bilinear estimate} 
For the broad part in (\ref{s2s3e10}), we are confronted with bilinear terms with several uncontrollable transversality parameters $\mu_1$,...,$\mu_t$. To obtain desired bounds, we introduce the following bilinear restriction estimate independent of the transversality parameters. The proof relies on a classical $L^4$ argument, and has been derived by Guo--Oh \cite{GO22} for the special case $\mathbf{Q}=(\xi_1^2,\xi_2^2+\xi_1\xi_3)$.

\begin{theorem}\label{s3thm1}
    Let $\mathbf{Q}$ be an $n$-tuple of quadratic forms on $\mathbb{R}^d$ with $d\geq n$. Suppose that there exist several parameters $i_1,$...,$i_{d-n}$ such that
    \begin{equation}\label{s3the1e1}
        |J_{\mathbf{Q}}(\xi;i_1,...,i_{d-n})|\sim \prod_{j=1}^t |\xi_j|^{s_j},
    \end{equation}
    where $t\leq d$, $ s_j \in \mathbb{N}^+$ {\rm(}$1\leq j\leq t${\rm)}, and $\sum_{j=1}^t s_j=n$. Then we have 
    \begin{equation*}
        \sup_{1<\mu_1,...,\mu_t\leq R^{1/2}}BD_p(\mu_1,...,\mu_t,1,...,1;R)\lesssim 1
    \end{equation*}
    for any $p\geq \max_{j} s_j+3.$
\end{theorem}

\noindent \textit{Proof.}  By Definition \ref{s2def1}, we need to consider the bilinear restriction estimate
$$\Big\||E^{\mathbf{Q}}f_1 E^{\mathbf{Q}}f_2|^{\frac{1}{2}}\Big\|_{L^p(\mathbb{R}^{d+n})} \lesssim  \|f_1\|_{L^p}^{\frac{1}{2}} \|f_2\|_{L^p}^{\frac{1}{2}},$$
where $f_1$ and $f_2$ are supported in two separated rectangular boxes of dimensions $\mu_1^{-1}\times \cdots \times \mu_t^{-1}\times 1\times \cdots \times 1$.
When $p=\infty$, by H\"older's inequality, we trivially have
$$   \Big\||E^{\mathbf{Q}}f_1 E^{\mathbf{Q}}f_2|^{\frac{1}{2}}\Big\|_{L^\infty(\mathbb{R}^{d+n})} \leq \|f_1\|^{\frac{1}{2}}_{L^1}\|f_2\|^{\frac{1}{2}}_{L^1}\leq \mu_1^{-\frac{1}{2}}...\mu_t^{-\frac{1}{2}} \|f_1\|^{\frac{1}{2}}_{L^2}\|f_2\|^{\frac{1}{2}}_{L^2}.     $$
When $p=4$, we write
$$   E^{\mathbf{Q}}f_1(x) E^{\mathbf{Q}}f_2(x)= \int \int f_1(\xi) f_2(\xi') e\Big[x'(\xi+\xi')+x''(\mathbf{Q}(\xi)+\mathbf{Q}(\xi'))\Big] \dd\xi \dd\xi',   $$
with $x'=(x_1,...,x_d)$, $x''=(x_{d+1},...,x_{d+n})$, $\xi=(\xi_1,..,\xi_d)$ and $\xi'=(\xi_1',...,\xi_{d}')$. Apply the change of variables
\begin{equation}
    ~\begin{cases}
        ~	\eta'=\xi+\xi';   \\
        ~	\eta''=\mathbf{Q}(\xi)+\mathbf{Q}(\xi');  \\
        ~	\eta_{d+n+1}=\xi'_{i_1};  \\
        ~ \quad   \vdots    \\
        ~	\eta_{2d}=\xi'_{i_{d-n}},
    \end{cases}
\end{equation}
where $\eta'=(\eta_1,...,\eta_{d})$ and $\eta''=(\eta_{d+1},...,\eta_{d+n})$. Then we see that the Jacobian $J$ of the change of variables is just $J_{\mathbf{Q}}(\xi'-\xi;i_1,...,i_{d-n})$. Thus, we can compute the $L^4$ norm as follows:
\begin{align*}
    \Big\||E^{\mathbf{Q}}f_1 E^{\mathbf{Q}}f_2|^{\frac{1}{2}}\Big\|^4_{L^4(\mathbb{R}^{d+n})}&=\int  |E^{\mathbf{Q}}f_1 E^{\mathbf{Q}}f_2|^{2}   \\
    &= \int \left| \int \int g_1(\eta) g_2(\eta) e(x\cdot \tilde{\eta}) J^{-1}\dd\tilde{\eta}\dd\overline{\eta}  \right|^2 \dd x   \\
    & \leq \int \int |g_1|^2 |g_2|^2 J^{-2} \\
    & \leq \int \int |f_1|^2 |f_2|^2 J^{-1}  \\
    &  \lesssim \mu_1^{s_1}\cdots\mu_t^{s_t}  \|f_1\|^2_{L^2} \|f_2\|^2_{L^2},
\end{align*}
where $g_1$ and $g_2$ are appropriate functions, $\tilde{\eta}=(\eta_1,...,\eta_{d+n})$, $\overline{\eta}=(\eta_{d+n+1},,,.\eta_{2d})$, and $\eta = (\tilde{\eta},\overline{\eta})$. Here in the third line we first use H\"older's inequality in the variables $\overline{\eta}$ and then apply the Plancherel theorem in the variables $\tilde{\eta}$, in the fourth line we change variables back, and in the last line we apply the assumption (\ref{s3the1e1}). It follows that 
$$\Big\||E^{\mathbf{Q}}f_1 E^{\mathbf{Q}}f_2|^{\frac{1}{2}}\Big\|_{L^4(\mathbb{R}^{d+n})} \lesssim \mu_1^{\frac{s_1}{4}}\cdots\mu_t^{\frac{s_t}{4}}  \|f_1\|^{\frac{1}{2}}_{L^2} \|f_2\|^{\frac{1}{2}}_{L^2} . $$
Through H\"older's inequality, we can conclude that for $4\leq p \leq \infty$,
\begin{align*}
    \Big\||E^{\mathbf{Q}}f_1 E^{\mathbf{Q}}f_2|^{\frac{1}{2}}\Big\|_{L^p(\mathbb{R}^{d+n})} &\lesssim  \mu_1^{\frac{s_1+2}{p}-\frac{1}{2}}\cdots\mu_t^{\frac{s_t+2}{p}-\frac{1}{2}}  \|f_1\|^{\frac{1}{2}}_{L^2} \|f_2\|^{\frac{1}{2}}_{L^2} \\
    & \lesssim \mu_1^{\frac{s_1+3}{p}-1}\cdots\mu_t^{\frac{s_t+3}{p}-1}  \|f_1\|^{\frac{1}{2}}_{L^p} \|f_2\|^{\frac{1}{2}}_{L^p} .
\end{align*}
As long as we take $p\geq \max_{j} s_j+3$, by noting that $\mu_j > 1$ for all $1\leq j\leq t$, we have 
$$\Big\||E^{\mathbf{Q}}f_1 E^{\mathbf{Q}}f_2|^{\frac{1}{2}}\Big\|_{L^p(\mathbb{R}^{d+n})} \lesssim   \|f_1\|^{\frac{1}{2}}_{L^p} \|f_2\|^{\frac{1}{2}}_{L^p}$$
as desired. \qed

\subsection{Decoupling}

For the narrow part in (\ref{s2s3e10}), we plan to apply different techniques, such as the locally constant property, induction on scale, adaptive decoupling and so on. Here we state some classical decoupling results that will be needed in the next section.

\begin{lemma}
    [$\ell^p$ decoupling for hypersurfaces with nonzero Gaussian curvature, \cite{BD17}] 
    Let $F:\mathbb{R}^m \rightarrow \mathbb{C}$ be such that ${\rm supp}\widehat{F}\subset N_{R^{-1}}(S)$, where $S$ is a compact $C^2$ hypersurface in $\mathbb{R}^m$ with nonzero Gaussian curvature. Divide this neighborhood into slabs $\theta$ with $m-1$ long directions of length $R^{-1/2}$ and one short direction of length $R^{-1}$. Write $F=\sum_\theta F_\theta$, where $\widehat{F}_\theta=\widehat{F} \chi_\theta$. Then 
    \begin{equation}\label{d1}
        \|F\|_{L^p(\mathbb{R}^m)} \lesssim_\epsilon R^{\beta(p)+\epsilon} \Big( \sum_{\theta} \|F_\theta\|_{L^p(\mathbb{R}^m)}^p  \Big)^{\frac{1}{p}},
    \end{equation}
    where $\beta(p)=\frac{m-1}{4}-\frac{m-1}{2p}$ when $2\leq p < \frac{2(m+1)}{m-1}$ and $\beta(p)=\frac{m-1}{2}-\frac{m}{p}$ when $\frac{2(m+1)}{m-1}\leq p \leq \infty$.
\end{lemma}

For our purposes, we will actually use the localized version of (\ref{d1}), i.e.,
\begin{equation}\label{d10}
        \|F\|_{L^p(B_R^m)} \lesssim_\epsilon R^{\beta(p)+\epsilon} \Big( \sum_{\theta} \|F_\theta\|_{L^p(B_R^m)}^p  \Big)^{\frac{1}{p}}.
    \end{equation}
It is well-known that (\ref{d1}) and (\ref{d10}) are equivalent, and the reader may consult \cite[Proposition 9.15]{demeter2020} for more details. Strictly speaking, we should consider $\|F_\theta\|_{L^p(w_{B_R})}$ rather than $\|F_\theta\|_{L^p(B_R)}$ for some weighted function $w_{B_R}$ on the right side of the inequality in (\ref{d10}). We omit such technical details since corresponding treatments of the two cases are essentially the same.

\begin{lemma} [Flat decoupling, \cite{BD17}]
    Let $B$ be a rectangular box in $\mathbb{R}^m$, and $B_1,...,B_L$ be a partition of $B$ into congruent boxes that are translates of each other. Write $F=\sum_j F_{j}$, where $\widehat{F}_{j}=\widehat{F} \chi_{B_j}$. Then for any $2\leq p,q\leq \infty$, we have
    \begin{equation}\label{d2}
        \|F\|_{L^p(\mathbb{R}^m)} \lesssim L^{1-\frac{1}{p}-\frac{1}{q}}\Big( \sum_{j} \|F_j\|_{L^p(\mathbb{R}^m)}^q  \Big)^{\frac{1}{q}}.
    \end{equation}
\end{lemma}

\section{The proof of Theorem \ref{main thm}}\label{section4}

In this section, we start to prove Theorem \ref{main thm}. Suppose that $\mathbf{Q}$ is an $n$-tuple of quadratic forms on $\mathbb{R}^d$ as in Theorem \ref{main thm}. By the trivial $L^1 \rightarrow L^\infty $ estimate and interpolation, it suffices to show 
\begin{align}\label{eq:goal} \|E^{\mathbf{Q}}f\|_{L^p(\mathbb{R}^{d+n})} \lesssim \|f\|_{L^p(\mathbb{R}^d)},
\end{align}
for $p> \max_j w_j+2$. Fix any such $p$ from now on. From the condition that $J_{\mathbf{Q}}(\xi;i_1,...,i_{d-n}) \not\equiv 0$ for some $i_1,$...,$i_{d-n}$, we know that $\mathbf{Q}$ must be linearly independent. By \cite[Theorem 3.1]{CLMP25}, the associated surface $S_{\mathbf{Q}}$ has nontrivial uniform Fourier decay. Thus we can apply the epsilon removal lemma in \cite[Corollary 7.5]{CLMP25}, and (\ref{eq:goal}) can be further reduced to its localized version: For any $0<\epsilon\ll1$,
$$    \|E^{\mathbf{Q}}f\|_{L^p(B_R^{d+n})} \lesssim R^\epsilon\|f\|_{L^p(\mathbb{R}^d)}$$
holds for all $R\geq1$. Next, by Definition \ref{s2def1}, it is equivalent to show: For any $0<\epsilon\ll1$,
\begin{equation}\label{s4e1}
    D_p(R) \lesssim R^\epsilon
\end{equation}
holds for all $R\geq 1$. 
Note that by the assumption $\max_j w_j\neq 1$, we always have $p > 4$. Besides, we can assume that $R \gg 1$. Otherwise, we can take a large implicit constant in (\ref{s4e1}), and then (\ref{s4e1}) holds trivially. Recall that there are two cases in Theorem~\ref{main thm}, and let us treat them separately. 

\vskip0.3cm

(1) We are in the monomial case, that is,
$$\mathbf{Q}(\xi)=(\xi_{\lambda_1} \xi_{\lambda_2},\xi_{\lambda_3}\xi_{\lambda_4},...,\xi_{\lambda_{2n-1}}\xi_{\lambda_{2n}}).$$

\textbf{Claim}: For any $0<\epsilon < 1-\frac{\max_j w_j+2}{p}$ (recall that $p>\max_j w_j +2$), (\ref{s4e1}) holds for all $R\gg1$. 

By our assumption, there exist $i_1,$...,$i_{d-n}$ such that $J_{\mathbf{Q}}(\xi;i_1,...,i_{d-n}) \not\equiv 0$. So, by Theorem~\ref{thm:Jacobian} and the reduction explained at the beginning of Section~\ref{section3}, we can assume that \begin{equation}\label{s4c1e0}
    |J_{\mathbf{Q}}(\xi;i_1,...,i_{d-n})|\sim \prod_{j=1}^t |\xi_j|^{s_j}
\end{equation}
with $s_j \in \N^+$ {\rm(}$1\leq j\leq t${\rm)}, $\sum_{j=1}^t s_j=n$, and 
\begin{equation}\label{s4add1}
    \max_j s_j \leq \max_j w_j-1.
\end{equation}
By (\ref{s2s3e10}) (replace $\epsilon$ with $\epsilon/2$), we obtain
\begin{equation}\label{s4c1e1}
D_p(R) \lesssim R^{\frac{\epsilon}{2}}\sup_{1<\mu_1,...,\mu_t\leq R^{1/2}} BD_p(\mu_1,...,\mu_t,1,...,1;R)+  R^{\frac{\epsilon}{2}}\sum_{j=1}^t D_p(1,...,1,\underset{\substack{\uparrow \\ j}}{R^\frac{1}{2}},1,...,1;R).
\end{equation}

For the first bilinear term in (\ref{s4c1e1}), we can use Theorem \ref{s3thm1} to conclude that
$$R^{\frac{\epsilon}{2}}\sup_{1<\mu_1,...,\mu_t\leq R^{1/2}} BD_p(\mu_1,...,\mu_t,1,...,1;R) \ll R^\epsilon,$$
since $p > \max_j w_j +2 \geq \max_j s_j+3$ by (\ref{s4add1}).

For the remaining linear terms in (\ref{s4c1e1}), we take the $j$th term 
$$   R^{\frac{\epsilon}{2}} D_p(1,...,1,\underset{\substack{\uparrow \\ j}}{R^\frac{1}{2}},1,...,1;R) $$
as an example. Without loss of generality, we may assume that $\supp\,f$ in this setting is in
$$   [0,1]\times \cdots\times [0,1]\times [0,\underset{\substack{\uparrow \\ j}}R^{-\frac{1}{2}}]\times [0,1]\times \cdots \times[0,1].    $$
Since every component of $\mathbf{Q}$ is a quadratic monomial, the related surface $S_\mathbf{Q}$ always stays unchanged after rescaling. If we use the rescaling $\xi_j \rightarrow \xi_j/R^{1/2}$, then
$$ D_p(1,...,1,\underset{\substack{\uparrow \\ j}}{R^\frac{1}{2}},1,...,1;R) \leq R^{\frac{w_j+2}{2p}-\frac{1}{2}} D_p(R)\lesssim R^{\frac{w_j+2}{2p}-\frac{1}{2}+\epsilon},$$
where we cover $B_R$ with balls of scale $R/2$ and apply the induction hypothesis (\ref{s4e1}). Thus
\begin{equation*}
    R^{\frac{\epsilon}{2}} D_p(1,...,1,\underset{\substack{\uparrow \\ j}}{R^\frac{1}{2}},1,...,1;R) \lesssim R^{\frac{w_j+2}{2p}-\frac{1}{2}+\frac{3}{2}\epsilon}\ll R^{\epsilon},
\end{equation*}
due to the fact $0<\epsilon < 1-\frac{\max_j w_j+2}{p}  \leq 1-\frac{w_j+2}{p}$.

Combining the estimates for all terms in (\ref{s4c1e1}), we validate \textbf{Claim}, which completes the proof.

\vskip0.3cm

(2) In this case, we assume that 
$$\mathbf{Q}(\xi)=(\xi_{\lambda_1} \xi_{\lambda_2},\xi_{\lambda_3} \xi_{\lambda_4},...,\xi_{\lambda_{2n-3}}\xi_{\lambda_{2n-2}},\xi_{\lambda_{2n-1}}\xi_{\lambda_{2n}}+P(\xi))$$
is an  $n$-tuple of quadratic forms defined in $\mathbb{R}^d$. Since $P(\xi)\not\equiv0$ is independent of $\{ \xi_{\lambda_j}\}_{j=1}^{2n}$, by an invertible linear transformation, we can assume that
$$  \mathbf{Q}(\xi)=(\xi_{\lambda_1} \xi_{\lambda_2},\xi_{\lambda_3}\xi_{\lambda_4},...,\xi_{\lambda_{2n-3}}\xi_{\lambda_{2n-2}},\xi_{\lambda_{2n-1}}\xi_{\lambda_{2n}} \pm \xi_{d-k+1}^2\pm\cdots \pm \xi_d^2),      $$
for some $k\geq 1$. Here each $\xi_{\lambda_j}$ ($j=1,...,2n$) comes from $\{\xi_1,...,\xi_{d-k}\}$. Recall that $\tilde{\mathbf{Q}}(\xi)=(\xi_{\lambda_1} \xi_{\lambda_2},...,\xi_{\lambda_{2n-3}}\xi_{\lambda_{2n-2}},\xi_{\lambda_{2n-1}}\xi_{\lambda_{2n}})$ is an $n$-tuple of quadratic forms on $\mathbb{R}^{\tilde{d}}$, so $d=\tilde{d}+k $.

\textbf{Claim}: There exists $\epsilon_0>0$, such that for any $0<\epsilon < \epsilon_0$, (\ref{s4e1}) holds for all $R\gg1$.

The admissible value of $\epsilon_0$ will be determined along the way of our proof.

First, we can assume that $\xi_{\lambda_{2n-1}} $ or $\xi_{\lambda_{2n}}$ must appear in the variables of 
$$\mathbf{Q}'(\xi)=(\xi_{\lambda_1} \xi_{\lambda_2},\xi_{\lambda_3} \xi_{\lambda_4},...,\xi_{\lambda_{2n-3}}\xi_{\lambda_{2n-2}}).$$
Otherwise, we can decompose $\mathbf{Q}$ into two independent quadratic forms $\mathbf{Q}=(\mathbf{Q}',\mathbf{Q}'')$ with
$$   \mathbf{Q}''(\xi)=(\xi_{\lambda_{2n-1}}\xi_{\lambda_{2n}}  \pm \xi_{d-k+1}^2\pm\cdots \pm \xi_d^2).$$ 
Note that $\mathbf{Q}'$ has dimension $d'\geq \tilde{d}-2$ and codimension $n'=n-1$. The condition $\tilde{d}\geq n$ itself does not guarantee $d'\geq n'$. However, by the Jacobian condition for $\tilde{\mathbf{Q}}$  and Theorem~\ref{thm:inherit}, we indeed have $d' \geq n'$. Moreover, $J_{\mathbf{Q}'}(\xi;i_1,...,i_{d'-n'}) \not\equiv 0$ for some $i_1,$...,$i_{d'-n'}$. Denote the number of occurrences of $\xi_j$ in $\{ \xi_{\lambda_j}:j=1,...,2n-2  \}$ by $w_j'$. Now we begin to show that \begin{equation}\label{eq:Q'}
	\|E^{\mathbf{Q}'}g\|_{L^p(B_R^{d'+n'})}\lesssim R^\epsilon \|g\|_{L^p(\mathbb{R}^{d'})}
\end{equation}
for any $g\in L^p(\mathbb{R}^{d'})$. There are two possibilities. If $\max_j w_j'\not=1$, then we can directly apply the result of Case (1) in Theorem \ref{main thm} to $\mathbf{Q}'$ to obtain (\ref{eq:Q'}) due to $p>\max_j w_j+2$. On the other hand, if $\max_j w_j'=1$, then $\mathbf{Q}'$ must be a tensor product of $n'$ hyperbolic paraboloids. Without loss of generality, we can assume 
$$   \mathbf{Q}'(\xi)=(\xi_1\xi_2,\xi_3\xi_4,...,\xi_{2n'-1}\xi_{2n'}).  $$
Using the restriction estimate for the hyperbolic paraboloid in $\mathbb{R}^3$ (see \cite{GHI19}) for each monomial component, we can still obtain (\ref{eq:Q'}) due to $p>\max_j w_j+2\geq 4$. Thus, (\ref{eq:Q'}) always holds true. Then by applying restriction estimates for the hyperbolic paraboloids to $\mathbf{Q}''$ (note $p>4$, see \cite{GHI19}), we immediately get
$$  \|E^{\mathbf{Q}}f\|_{L^p(B_R)} =\|E^{\mathbf{Q}'}(E^{\mathbf{Q}''}f)\|_{L^p(B_R)} \lesssim R^\epsilon \|f\|_{L^p}, $$
which validates \textbf{Claim}. So, we can indeed assume that $\xi_{\lambda_{2n-1}} $ or $\xi_{\lambda_{2n}}$ must appear in 
$\mathbf{Q}'$.

Now we start to perform the broad--narrow analysis. By the Jacobian condition for $\tilde{\mathbf{Q}}$ and Theorem~\ref{thm:Jacobian}, we can assume that 
\begin{equation*}
	|J_{\tilde{\mathbf{Q}}}(\xi;i_1,...,i_{\tilde{d}-n})|\sim \prod_{j=1}^t |\xi_j|^{s_j}
\end{equation*}
with $s_j \in \N^+$ {\rm(}$1\leq j\leq t${\rm)}, $\sum_{j=1}^t s_j=n$, and 
\begin{equation}\label{s4add2}
	\max_j s_j\leq    \max_{j} w_j-1.  
\end{equation}
Via a simple calculation, we have 
$$ |J_{\mathbf{Q}}(\xi;i_1,...,i_{\tilde{d}-n},\tilde{d}+1,\tilde{d}+2,...,\tilde{d}+k)|\sim |J_{\tilde{\mathbf{Q}}}(\xi;i_1,...,i_{\tilde{d}-n})|\sim  |\xi_{1}|^{s_{1}}\cdots |\xi_{t}|^{s_{t}}.   $$
So, we can use (\ref{s2s3e10}) (replace $\epsilon$ with $\epsilon/2$) to obtain
\begin{align}
	D_p(R) \lesssim  R^{\frac{\epsilon}{2}}\sup_{1<\mu_1,...,\mu_t\leq R^{1/2}} BD_p(\mu_1,...,\mu_t,1,...,1;R)  
	+R^{\frac{\epsilon}{2}}\sum_{j=1}^t D_p(1,...,1,\underset{\substack{\uparrow \\ j}}{R^\frac{1}{2}},1,...,1;R),  \label{s4c2e1}
\end{align}
with $t\leq \tilde{d}=d-k$. For the first bilinear term in (\ref{s4c2e1}), we can use Theorem \ref{s3thm1} to conclude that
$$R^{\frac{\epsilon}{2}}\sup_{1<\mu_1,...,\mu_t\leq R^{1/2}} BD_p(\mu_1,...,\mu_t,1,...,1;R) \ll R^\epsilon,$$
since $p > \max_j w_j +2 \geq \max_j s_j+3$ by (\ref{s4add2}). For the remaining linear terms in (\ref{s4c2e1}), we take the $j$th term 
$$   R^{\frac{\epsilon}{2}}  D_p(1,...,1,\underset{\substack{\uparrow \\ j}}{R^\frac{1}{2}},1,...,1;R) $$
as an example. As in the proof for Case (1), we may assume that $\supp\,f$ in this setting is in 
$$   [0,1]\times \cdots\times [0,1]\times [0,\underset{\substack{\uparrow \\ j}}R^{-\frac{1}{2}}]\times [0,1]\times \cdots \times[0,1].    $$

We first suppose that $\xi_{\lambda_{2n-1}}\xi_{\lambda_{2n}}$ is a mixed term. Then, without loss of generality, we can assume that $\xi_{\lambda_{2n-1}}\xi_{\lambda_{2n}}=\xi_{d-k-1}\xi_{d-k}$, that is,
$$  \mathbf{Q}(\xi)=(\xi_{\lambda_1} \xi_{\lambda_2},\xi_{\lambda_3}\xi_{\lambda_4},...,\xi_{\lambda_{2n-3}}\xi_{\lambda_{2n-2}},\xi_{d-k-1}\xi_{d-k} \pm \xi_{d-k+1}^2\pm\cdots \pm \xi_d^2).      $$ 
There are two subcases for the exact position of $j$, and now we treat them in order.

 (i) If $1\leq j\leq d-k-2$, then $\xi_j$ does not appear in the last component of $\mathbf{Q}$, and each previous component is a quadratic monomial. So, we can use rescaling and induction on scale as in Case (1), and taking $\epsilon_0 < 1-\frac{\max_j w_j+2}{p}$ will close the induction.

 (ii) If $d-k-1\leq j\leq d-k$, then without loss of generality, we can assume $j=d-k-1$. 
 
 Let us suppose $k>1$ first. By the locally constant property, we have
 $$\|E^{\mathbf{Q}}f\|_{L^p(B_{R^{1/2}})}   \sim \|E^{\mathbf{Q}_1}f\|_{L^p(B_{R^{1/2}})},$$
 where
 $$  \mathbf{Q}_1(\xi)=(\xi_{\lambda_1} \xi_{\lambda_2},\xi_{\lambda_3}\xi_{\lambda_4},...,\xi_{\lambda_{2n-3}}\xi_{\lambda_{2n-2}},\pm \xi_{d-k+1}^2\pm\cdots \pm \xi_d^2).    $$
In fact, recall
\begin{align*}
  E^{\mathbf{Q}}f(x)=\int_{[0,1]^d} e\Big[x_1\xi_1+\cdots+x_{d}\xi_d+&x_{d+1}\xi_{\lambda_1} \xi_{\lambda_2}+\cdots+ x_{d+n-1}\xi_{\lambda_{2n-3}}\xi_{\lambda_{2n-2}}+  \\
  &  x_{d+n}\bigl(\xi_{d-k-1}\xi_{d-k} \pm \xi_{d-k+1}^2\pm\cdots \pm \xi_d^2\bigr)\Big] f(\xi)d\xi,  
\end{align*}
where $e(y):=e^{2\pi i y}$ for any $y\in \mathbb{R}$. By translation, without loss of generality, we can assume that $B_{R^{1/2}}$ is centered at $0$. Then by the ranges of $x$ and $\xi$: $x\in B(0,R^{1/2})$, $\xi_{d-k-1}\in [0,R^{-1/2}]$ and $\xi_{d-k}\in [0,1]$, we obtain $|x_{d+n}\xi_{d-k-1}\xi_{d-k}|\leq 1$. If we further shrink the region, for instance, by taking $\xi_{d-k}\in [0,1/10]$ (such treatment will only produce $O(1)$-loss in the implicit constant), then the oscillatory term $e(x_{d+n}\xi_{d-k-1}\xi_{d-k})\sim 1$, which can be neglected. We refer the reader to \cite[Lemma 6.1]{GWZ20} for the rigorous statement and proof of the locally constant property.

Since each $\xi_{\lambda_j}$ ($j=1,...,2n-2$) comes from $\{\xi_1,...,\xi_{d-k}\}$, we can decompose $\mathbf{Q}_1$ into two independent quadratic forms $\mathbf{Q}_1=(\mathbf{Q}',\mathbf{Q}''_1)$ with
$$     \mathbf{Q}_1''(\xi)= (\pm\xi_{d-k+1}^2\pm\cdots \pm \xi_d^2).    $$
Note that $E^{\mathbf{Q}_1}f=E^{\mathbf{Q}_1''}(E^{\mathbf{Q}'}f)$, and $\mathbf{Q}_1''$ denotes a hyperbolic paraboloid, we can apply (\ref{d10}) to $\mathbf{Q}_1''$ to obtain
 \begin{equation*}
  \|E^{\mathbf{Q}}f\|_{L^p(B_{R^{1/2}})} \lesssim R^{\frac{k}{4}-\frac{k+1}{2p}+\epsilon} \Big(\sum_{\theta} \|E^{\mathbf{Q}}f_\theta\|_{L^p(B_{R^{1/2}})}^p  \Big)^{\frac{1}{p}},   
 \end{equation*}
 where each $f_\theta$ is supported in a rectangular box $\theta$ of size 
 $$1\times \cdots\times 1\times \underset{\substack{\uparrow \\ d-k-1}}{R^{-\frac{1}{2}}}\times 1\times \underset{\substack{\uparrow \\ d-k+1}}{R^{-\frac{1}{4}}}\times \cdots \times R^{-\frac{1}{4}}.$$
 The constraints $p>4$ and $k>1$ guarantee that we are always in the supercritical range. 
Summing over all $B_{R^{1/2}}$ in $B_R$, we get 
\begin{equation}\label{s4c2e2}
   \|E^{\mathbf{Q}}f\|_{L^p(B_R)} \lesssim R^{\frac{k}{4}-\frac{k+1}{2p}+\epsilon} \Big(\sum_{\theta} \|E^{\mathbf{Q}}f_\theta\|_{L^p(B_R)}^p  \Big)^{\frac{1}{p}}.  
\end{equation}
By Definition \ref{s2def1}, this immediately implies
$$   D_p(1,...,1,\underset{\substack{\uparrow \\ d-k-1}}{R^\frac{1}{2}},1,...,1;R) \lesssim  R^{\frac{k}{4}-\frac{k+1}{2p}+\epsilon}  D_p(1,...,1,\underset{\substack{\uparrow \\ d-k-1}}{R^\frac{1}{2}},1,\underset{\substack{\uparrow \\ d-k+1}}{R^\frac{1}{4}},...,R^{\frac{1}{4}};R) .  $$

Next we use the change of variables
\begin{equation*}
    \xi_{d-k-1}\rightarrow\frac{\xi_{d-k-1}}{R^{1/2}},\quad \xi_{d-k+1}\rightarrow\frac{\xi_{d-k+1}}{R^{1/4}} , \quad... \quad,  \quad  \xi_{d} \rightarrow \frac{\xi_d}{R^{1/4}}.
\end{equation*}
Observe that such a rescaling leaves the surface $S_{\mathbf{Q}}$ unchanged. Thus, we get
$$  D_p(1,...,1,\underset{\substack{\uparrow \\ d-k-1}}{R^\frac{1}{2}},1,\underset{\substack{\uparrow \\ d-k+1}}{R^\frac{1}{4}},...,R^{\frac{1}{4}};R) \lesssim R^{\frac{w_{d-k-1}+k+2}{2p}-\frac{k+2}{4}}D_p(R)\lesssim R^{\frac{w_{d-k-1}+k+2}{2p}-\frac{k+2}{4}+\epsilon},    $$
where we cover $B_R$ with balls of scale $R/2$ and apply the induction hypothesis (\ref{s4e1}). Combining all estimates above, one concludes
$$ R^{\frac{\epsilon}{2}} D_p(1,...,1,\underset{\substack{\uparrow \\ d-k-1}}{R^\frac{1}{2}},1,...,1;R) \lesssim  R^{\frac{w_{d-k-1}+1}{2p}-\frac{1}{2}+\frac{5}{2}\epsilon} \ll R^\epsilon,           $$
as long as we take $\epsilon_0 < \frac{1}{3}\Big(1 - \frac{w_{d-k-1} + 1}{p}\Big)$. 
 
 If $k=1$, then still by (\ref{d10}), our (\ref{s4c2e2}) should be replaced by 
 \begin{equation*}
  \|E^{\mathbf{Q}}f\|_{L^p(B_R)} \lesssim R^{\beta(p)+\epsilon} \Big(\sum_{\theta} \|E^{\mathbf{Q}}f_\theta\|_{L^p(B_R)}^p  \Big)^{\frac{1}{p}},   
 \end{equation*}
 where $\beta(p)=\frac{1}{8}-\frac{1}{4p}$ when $4\leq p<6$ and $\beta(p)=\frac{1}{4}-\frac{1}{p}$ when $p\geq 6$. Same arguments as in the case $k>1$ can also lead to desired bounds by taking $\epsilon_0$ suitably small.

Now that we have validated \textbf{Claim} when $\xi_{\lambda_{2n-1}}\xi_{\lambda_{2n}}$ is a mixed term, it remains to consider the case when $\xi_{\lambda_{2n-1}}\xi_{\lambda_{2n}}$ is a square term. In this case, without loss of generality, we can assume that $\xi_{\lambda_{2n-1}}\xi_{\lambda_{2n}}=\xi^2_{d-k}$, that is,
$$  \mathbf{Q}(\xi)=(\xi_{\lambda_1} \xi_{\lambda_2},\xi_{\lambda_3}\xi_{\lambda_4},...,\xi_{\lambda_{2n-3}}\xi_{\lambda_{2n-2}},\xi_{d-k}^2 \pm \xi_{d-k+1}^2\pm\cdots \pm \xi_d^2).      $$ 
There are also two subcases for the exact position of $j$:
 
(i) If $1\leq j\leq d-k-1$, we can repeat the arguments in Case (1) as before to close the induction.
 
 (ii) If $ j= d-k$, then by the locally constant property, the original form $\mathbf{Q}$ can be reduced to
 $$  \mathbf{Q}_2(\xi)=(\xi_{\lambda_1} \xi_{\lambda_2},\xi_{\lambda_3}\xi_{\lambda_4},...,\xi_{\lambda_{2n-3}}\xi_{\lambda_{2n-2}}, \xi_{d-k+1}^2\pm\cdots \pm \xi_{d}^2).     $$
 Since each $\xi_{\lambda_j}$ ($j=1,...,2n-2$) comes from $\{\xi_1,...,\xi_{d-k}\}$, we can decompose $\mathbf{Q}_2$ into two independent quadratic forms $\mathbf{Q}_2=(\mathbf{Q}',\mathbf{Q}''_2)$ with
 $$     \mathbf{Q}_2''(\xi)= (\xi_{d-k+1}^2\pm\cdots \pm \xi_d^2).    $$
 In view of the initial reduction ($\xi_{\lambda_{2n-1}} $ or $\xi_{\lambda_{2n}}$ must appear in 
$\mathbf{Q}'$), we can assume that the variable $\xi_{d-k}$ appears in $\mathbf{Q}'$. Therefore, $\mathbf{Q}'$ has dimension $d'=\tilde{d}$ and codimension $n'=n-1$. From $\tilde{d}\geq n\geq 2$, we get $d'> n'\geq 1$. Besides, by the Jacobian condition for $\tilde{\mathbf{Q}}$  and Theorem~\ref{thm:inherit}, one concludes that $J_{\mathbf{Q}'}(\xi;i_1,...,i_{d'-n'}) \not\equiv 0$ for some $i_1,$...,$i_{d'-n'}$. Via the same argument as in the proof of (\ref{eq:Q'}), we always have 
 \begin{equation*}
 	\|E^{\mathbf{Q}'}g\|_{L^p(B_R^{d'+n'})}\lesssim R^\epsilon \|g\|_{L^p(\mathbb{R}^{d'})}
 \end{equation*}
 for any $g\in L^p(\mathbb{R}^{d'})$ due to $p>\max_j w_j+2$. Then by applying restriction estimates for the hyperbolic paraboloids to $\mathbf{Q}_2''$ (note $p>4$, see \cite{GHI19}), we get
 $$  \|E^{\mathbf{Q}}f\|_{L^p(B_R)} \sim \|E^{\mathbf{Q}_2}f\|_{L^p(B_R)}  =\|E^{\mathbf{Q}'}(E^{\mathbf{Q}_2''}f)\|_{L^p(B_R)} \lesssim R^\epsilon \|f\|_{L^p}. $$
Therefore, \textbf{Claim} is validated in all cases, which completes the proof. \qed

\vskip0.3cm

So far, we have finished the proof of the estimate (\ref{main thm eq2}). And we still need to show that it is sharp up to the endpoint. The argument is essentially the same as that in \cite[section 3]{GO22}, so we only sketch it here. We first consider Case (1) of Theorem~\ref{main thm}, that is,
$$\mathbf{Q}(\xi)=(\xi_{\lambda_1} \xi_{\lambda_2},\xi_{\lambda_3}\xi_{\lambda_4},...,\xi_{\lambda_{2n-1}}\xi_{\lambda_{2n}}).$$
Take the characteristic function $\chi_{[0,R^{-t_1}]\times \cdots\times [0,R^{-t_d}]}$ with $0\leq t_1,...,t_d\leq 1$. By the locally constant property, we have 
\begin{equation}\label{nese1}
    \Big|E^\mathbf{Q}(\chi_{[0,R^{-t_1}]\times \cdots\times [0,R^{-t_d}]})(x)\Big|\sim R^{-t_1-\cdots-t_d},\quad \quad x\in T, 
\end{equation}
where 
\begin{align*}
    T:=\Big\{ x:|x_1|\leq \frac{1}{100}R^{t_1},...,&|x_d| \leq \frac{1}{100}R^{t_d},\\
    &|x_{d+1}|\leq \frac{1}{100}R^{t_{\lambda_1}+t_{\lambda_2}} ,..., |x_{d+n}|\leq \frac{1}{100}R^{t_{\lambda_{2n-1}}+t_{\lambda_{2n}}}\Big\}.
\end{align*}
By the definition of $w_j$, we note
\begin{equation}\label{nese2}
    |T|\sim R^{t_1+\cdots +t_d+t_{\lambda_1}+\cdots+t_{\lambda_{2n}}}=R^{\sum_{j=1}^d (w_j+1)t_j}.
\end{equation}
Divide $[0,1]^d$ into $\{\tau\}$, where each $\tau=I_1\times\cdots\times I_d$ is a rectangular box with $|I_1|=R^{-t_1},...,|I_d|=R^{-t_d}$. Let $\epsilon=(\epsilon_\tau)$ be a sequence of independent random variables where $\epsilon_\tau$ takes values $\pm 1$ with equal probability. Define 
$$  f=\sum_{\tau} \epsilon_\tau \chi_\tau , $$
and then (\ref{main aim}) becomes
$$  \Big\|E^{\mathbf{Q}} (\sum_{\tau} \epsilon_\tau \chi_\tau)   \Big\|_{L^q} \lesssim 1.  $$
Using Khintchine's inequality and Holder's inequality, we get 
$$   1 \gtrsim  \mathbb{E} \Big\|E^{\mathbf{Q}} (\sum_{\tau} \epsilon_\tau \chi_\tau)   \Big\|^q_{L^q}\sim \int \Big(\sum_\tau |E^\mathbf{Q}\chi_\tau|^{2}\Big)^{\frac{q}{2}}\geq \int \sum_\tau |E^\mathbf{Q}\chi_\tau|^{q} .  $$
It follows from (\ref{nese1}) that
$$   1 \gtrsim   R^{\sum_{j=1}^d (w_j+2)t_j} R^{-(\sum_{j=1}^d t_j)q}.  $$
And so
$$    q\geq \frac{\sum_{j=1}^d w_jt_{j}}{\sum_{j=1}^d t_j}+2. $$
On the other hand, we can also test (\ref{main aim}) on a single $f_\tau$: 
$$ R^{-(\sum_{j=1}^d t_j)+\frac{\sum_{j=1}^d (w_j+1)t_j}{q}} \lesssim \|E^{\mathbf{Q}}f_\tau\|_{L^q}\lesssim \|f_\tau\|_{L^p}\sim R^{-\frac{\sum_{j=1}^d t_j}{p}},    $$
which leads to the constraint
$$  \frac{1}{p}+\frac{\sum_{j=1}^d (w_j+1)t_j}{\sum_{j=1}^d t_j}\cdot\frac{1}{q} \leq 1.  $$
Thus, we obtain the necessary condition
\begin{equation}\label{s4ef}
  q\geq \frac{\sum_{j=1}^d w_jt_{j}}{\sum_{j=1}^d t_j}+2,\quad \frac{1}{p}+\frac{\sum_{j=1}^d (w_j+1)t_j}{\sum_{j=1}^d t_j}\cdot\frac{1}{q} \leq 1.  
\end{equation}
Without loss of generality, we assume that $\max_j w_j=w_1$. To obtain optimal necessary condition, we take $t_1=1$ and $t_j=0$ for $j\geq 2$, and then (\ref{s4ef}) just matches (\ref{main thm eq2}) up to the endpoint.

Next we consider Case (2) of Theorem~\ref{main thm}, that is,
$$\mathbf{Q}(\xi)=(\xi_{\lambda_1} \xi_{\lambda_2},\xi_{\lambda_3} \xi_{\lambda_4},...,\xi_{\lambda_{2n-3}}\xi_{\lambda_{2n-2}},\xi_{\lambda_{2n-1}}\xi_{\lambda_{2n}}+P(\xi)).$$
Since $P(\xi)\not\equiv0$ is independent of $\{ \xi_{\lambda_j}\}_{j=1}^{2n}$, by an invertible linear transformation, we can assume that
$$  \mathbf{Q}(\xi)=(\xi_{\lambda_1} \xi_{\lambda_2},\xi_{\lambda_3}\xi_{\lambda_4},...,\xi_{\lambda_{2n-3}}\xi_{\lambda_{2n-2}},\xi_{\lambda_{2n-1}}\xi_{\lambda_{2n}} \pm \xi_{d-k+1}^2\pm\cdots \pm \xi_d^2),      $$
for some $1\leq k\leq d-n$. If we test the characteristic function $\chi_{[0,R^{-t_1}]\times \cdots\times [0,R^{-t_d}]}$ with $0\leq t_1,...,t_d\leq 1$ again, then we can get
\begin{equation*}
    \Big|E^\mathbf{Q}(\chi_{[0,R^{-t_1}]\times \cdots\times [0,R^{-t_d}]})(x)\Big|\sim R^{-t_1-\cdots-t_d},\quad \quad x\in \tilde{T}, 
\end{equation*}
where 
\begin{align*}
    \tilde{T}:=\Big\{ x:|x_1|\leq \frac{1}{100}R^{t_1},...,&|x_d| \leq \frac{1}{100}R^{t_d},\\
    &|x_{d+1}|\leq \frac{1}{100}R^{t_{\lambda_1}+t_{\lambda_2}} ,..., |x_{d+n}|\leq \frac{1}{100}R^{       \min\{t_{\lambda_{2n-1}}+t_{\lambda_{2n}} ,2t_{d-k+1},...,2t_d        \}   }\Big\}.
\end{align*}
Without loss of generality, we may assume that $\max_j w_j=w_1$. The condition $\lambda_{2n-1},\lambda_{2n}\notin \{ k: \max_j w_j= w_k \}$ guarantees $\lambda_{2n-1}\not=1$ and $\lambda_{2n}\not=1$. When we take $t_1=1$ and $t_j=0$ for $j\geq 2$, there is
$$  |\tilde{T}|\sim R^{t_1+\cdots +t_d+t_{\lambda_1}+\cdots+t_{\lambda_{2n-2}}+  \min\{t_{\lambda_{2n-1}}+t_{\lambda_{2n}} ,2t_{d-k+1},...,2t_d        \}  }=R^{w_1+1},  $$
which exactly matches (\ref{nese2}). Repeating the same argument as in Case (1), we can get the same necessary condition (\ref{s4ef}).

\section{Final remarks and more applications}\label{section5}

In this final section, we give some additional remarks and applications regarding our method and Theorem \ref{main thm}.

First, by the method of this paper, we can recover essentially sharp restriction estimates for all quadratic forms with $d=n=2$. In what follows, we shall first give a complete classification for all quadratic forms with $d=n=2$. To achieve this purpose, we need to introduce a technical quantity $\mathfrak{d}_{d',n'}(\mathbf{Q})$. Guo--Oh--Zhang--Zorin-Kranich \cite{guo2023decoupling} used this quantity to measure ``curvature" properties of the quadratic surface $S_{\mathbf{Q}}$, and obtained sharp decoupling inequalities for all higher codimensional quadratic surfaces.

\begin{definition}[\cite{guo2023decoupling}]
For a tuple $\mathbf{Q}(\xi)=(Q_1(\xi),...,Q_n(\xi))$ of quadratic forms on $\mathbb{R}^d$, we define
$$  {\rm NV}(\mathbf{Q}) := \# \{    1 \leq d' \leq d : \partial_{\xi_{d'}} Q_{n'} \not\equiv 0 \text{~for~some~} 1\leq n'\leq n \}      .  $$
For $0 \leq d' \leq d$ and $0\leq n'\leq n$, we define
\begin{equation}\label{s2def2e1}
\mathfrak{d}_{d',n'}(\mathbf{Q}):= \inf_{\substack{   M \in \mathbb{R}^{d\times d} \\ {\rm rank}(M)=d'   }}  \inf_{\substack{   M' \in \mathbb{R}^{n\times n'} \\ {\rm rank}(M')=n'   }} {\rm NV}((\mathbf{Q}\circ M)\cdot M') ,
\end{equation}
where $\mathbf{Q}\circ M$ is the composition of $\mathbf{Q}$ and $M$. We say that $\mathbf{Q}$ and $\mathbf{Q}'$ are equivalent and write $\mathbf{Q}\equiv\mathbf{Q}'$ if there exist two invertible real matrices $M_1 \in \mathbb{R}^{d\times d}$ and $M_2 \in \mathbb{R}^{n\times n}$ such that 
$$  \mathbf{Q}'(\xi) =  \mathbf{Q}(M_1 \cdot \xi) \cdot M_2, \quad \quad \forall~ \xi \in \mathbb{R}^d. $$   
\end{definition}
\begin{theorem}
Let $\mathbf{Q}=(Q_1,Q_2)$ be a $2$-tuple of quadratic forms on $\mathbb{R}^2$. Suppose that $\mathbf{Q}$ does not miss any variable and linearly independent, then
\begin{itemize}
\item[${\rm(a)}$] if $\mathbf{Q}$ satisfies $\mathfrak{d}_{2,1}(\mathbf{Q})=1$ and $\mathfrak{d}_{1,2}(\mathbf{Q})=0$, then $\mathbf{Q}(\xi)\equiv (\xi_1^2,\xi_1\xi_2)$, and the estimate ${\rm(\ref{main aim})}$ holds for
\begin{equation}\label{s6e1}
 q > 5,\quad \frac{1}{p}+\frac{ 4}{q}<1;
\end{equation}
\item[${\rm(b)}$] if $\mathbf{Q}$ satisfies $\mathfrak{d}_{2,1}(\mathbf{Q})=1$ and $\mathfrak{d}_{1,2}(\mathbf{Q})=1$, then $\mathbf{Q}(\xi)\equiv (\xi_1^2,\xi_2^2)$, and the estimate ${\rm(\ref{main aim})}$ holds for
\begin{equation}\label{s6e2}
  q > 4,\quad \frac{1}{p}+\frac{3}{q}<1;
\end{equation}
\item[${\rm(c)}$] if $\mathbf{Q}$ satisfies $\mathfrak{d}_{2,1}(\mathbf{Q})=2$, then $\mathbf{Q}(\xi)\equiv (\xi_1 \xi_2,\xi_1^2-\xi_2^2)$, and the estimate ${\rm(\ref{main aim})}$ holds for
\begin{equation}\label{s6e3}
  q > 4,\quad \frac{1}{p}+\frac{3}{q}<1.
\end{equation}
\end{itemize}
Moreover, the estimates all above are sharp up to the endpoints.
\end{theorem}

\noindent \textit{Proof.} Suppose that $\mathbf{Q}$ satisfies $\mathfrak{d}_{2,1}(\mathbf{Q})=1$. By \cite[Lemma 2.2]{GORYZK19}, we have
\begin{equation}\label{00000}
    \mathfrak{d}_{2,1}(\mathbf{Q})=1 \Longleftrightarrow ~ \mathbf{Q}(\xi)\equiv(\xi_1^2,Q_2(\xi)).
\end{equation}
If $\mathfrak{d}_{1,2}(\mathbf{Q})=0$, by \cite[(1.14)]{guo2023decoupling}, this quantity must attain its minimum when $\xi_1=0$. We can assume that
$$   \mathbf{Q}(\xi)\equiv(\xi_1^2,\xi_1 L(\xi)),   $$
with a linear form $L$. Also note $\mathfrak{d}_{2,2}(\mathbf{Q})=2$, $L(\xi)$ must depend on the variable $\xi_2$. By a linear transformation, we get $\mathbf{Q}(\xi)\equiv (\xi_1^2,\xi_1\xi_2)$. If  $\mathfrak{d}_{1,2}(\mathbf{Q})=1$, we write
$$   \mathbf{Q}(\xi)\equiv(\xi_1^2,a\xi_2^2+b\xi_1\xi_2) ,  $$
for some $a$ and $b$. We know $a\neq 0$, otherwise it is just Case (a). If $b=0$, we have proved Case (b). If not, we add a multiple of $\xi_1^2$ to $Q_2$ such that the second term can form a perfect square. Using a linear transformation, we get $ \mathbf{Q}(\xi)\equiv(\xi_1^2,\xi_2^2).$ 

Suppose that $\mathbf{Q}$ satisfies $\mathfrak{d}_{2,1}(\mathbf{Q})=2$. We pick $M\in \mathbb{R}^{2\times 2}$ and $M'\in \mathbb{R}^{2\times 1}$ such that the equality in (\ref{s2def2e1}) is achieved with $d'=2,n'=1$. Applying suitable linear transformations, we can assume that $M=I_{2\times 2}$ and $M'=(1,0)^T$. Then $\mathfrak{d}_{2,1}(\mathbf{Q})=2$ implies that $Q_1$ depends on 2 variables. We use linear transformations again to diagonalize $Q_1$, then
$$   \mathbf{Q}(\xi)\equiv(\xi_1^2\pm\xi_2^2,Q_2(\xi)) .  $$
We first assume that $Q_1(\xi)=\xi_1^2+\xi_2^2$. Then we can reduce it to 
$$   \mathbf{Q}(\xi)\equiv(\xi_1^2+\xi_2^2,\lambda\xi_1\xi_2+\mu \xi_2^2) .  $$
If $\lambda=0$, it can be further reduced to $(\xi_1^2,\xi_2^2)$, which is just Case (b). If $\mu=0$, by adding a multiple of $\xi_1\xi_2$ to $Q_1$ such that $Q_1$ forms a perfect square, which is a contradiction to (\ref{00000}). Therefore, we can assume that $\lambda\neq 0$ and $\mu\neq 0$. In this case, we can add a multiple of $Q_1$ to $Q_2$ such that $Q_2$ forms a perfect square, which is a contradiction to (\ref{00000}) again. Next we assume that $Q_1(\xi)=\xi_1^2-\xi_2^2$. By a linear transformation, we can write
$$   \mathbf{Q}(\xi)\equiv(\xi_1\xi_2,a\xi_1^2+b\xi_2^2) .  $$
We can assume $a\neq 0$ and $b\neq 0$, otherwise it can be further reduced to Case (a). So
$$   \mathbf{Q}(\xi)\equiv(\xi_1\xi_2,\xi_1^2+b\xi_2^2) .  $$
If $b>0$, we add a multiple of $\xi_1\xi_2$ to $Q_2$ such that $Q_2$ can form a perfect square, which is a contradiction to (\ref{00000}). Therefore, $b<0$, and $\mathbf{Q}(\xi)\equiv (\xi_1 \xi_2,\xi_1^2-\xi_2^2)$. 

Now we show (\ref{s6e1})--(\ref{s6e3}). For (\ref{s6e1}) and (\ref{s6e2}), we have proved them in Case (1) in Theorem \ref{main thm}. For (\ref{s6e3}), though $ (\xi_1 \xi_2,\xi_1^2-\xi_2^2)$ is not included in Theorem \ref{main thm}, by noting that 
$$|J_{\mathbf{Q}}(\xi)|\sim \xi_1^2+\xi_2^2\geq \xi_1^2,$$
we can still apply Theorem \ref{s3thm1} and the arguments in Section \ref{section4}. The sharpness of these estimates can be shown by testing on characteristic functions as in the final part of the previous section. \qed

\begin{remark}
    As a historical remark, Guo--Oh {\rm\cite{GO22}} proved {\rm(\ref{s6e1})} for Case {\rm(a)}, and Christ {\rm\cite{Chr85}}  obtained endpoint results {\rm(}$q>4$, $\frac{1}{p}+\frac{3}{q} = 1${\rm)} for Case {\rm(b)} and Case {\rm(c)} via a different method. 
\end{remark}

Second, our method can also play a role in the nondegenerate case. We will take a typical class of examples to illustrate this point. Before this, we introduce a classical nondegenerate condition in the higher codimensional setting.
\begin{definition}[\cite{Chr82,Chr85,Moc96}]
For a tuple $\mathbf{Q}(\xi)=(Q_1(\xi),...,Q_n(\xi))$ of quadratic forms on $\mathbb{R}^d$, we say that $\mathbf{Q}$ satisfies the ${\rm(CM)}$ condition if
\begin{equation*}
\int_{\mathbb{S}^{n-1}} |\mathrm{det}(y_1 Q_1+...+y_n Q_n)|^{-\gamma}d\sigma(y) <\infty, \quad \quad \forall~0<\gamma<\frac{n}{d},
\end{equation*} 
where each $Q_j$ {\rm(}$j=1,...,n${\rm)} denotes the Hessian matrix associated with the quadratic form $Q_j(\xi)$, and  $d\sigma$ is the surface measure on the unit sphere $\mathbb{S}^{n-1}$.
\end{definition}

\begin{theorem}
Suppose that $\mathbf{Q}=(Q_1,Q_2)$ is a $2$-tuple of quadratic forms on $\mathbb{R}^3$ and satisfies the ${\rm(CM)}$ condition. Then we have
\begin{equation*}
  \|E^{\mathbf{Q}}f\|_{L^p(\mathbb{R}^5)} \lesssim  \|f\|_{L^p(\mathbb{R}^3)},
\end{equation*}
for $p> 4$.  
\end{theorem}

\noindent \textit{Proof.} When $d=3$ and $n=2$, Guo--Oh essentially proved that if $\textbf{Q}$ satisfies $\mathfrak{d}_{3,2}(\textbf{Q})=3$, then
\begin{equation*}
 \textbf{Q}  \text{~satisfies~the~} ({\rm CM}) \text{~condition~}\Longleftrightarrow ~ \mathfrak{d}_{3,1}(\textbf{Q})=2,~ \mathfrak{d}_{2,2}(\textbf{Q})=2.    
\end{equation*}
Readers may also consult \cite[Theorem 1.1]{CMW24} for this equivalence. Then via (c5) of \cite[Lemma 6.5]{CLMP25}, we know that there are only three cases:
$$ \mathbf{Q}(\xi)\equiv(\xi_1^2+\xi_2^2,\xi_2^2+\xi_3^2) \quad\text{or}\quad(\xi_1^2-\xi_2^2,\xi_2^2-\xi_3^2)\quad\text{or}\quad (\xi_1^2-\xi_2^2,\xi_1\xi_2+\xi_3^2). $$
We treat the first case; the other two are similar.

For $\mathbf{Q}(\xi)=(\xi_1^2+\xi_2^2,\xi_2^2+\xi_3^2)$, by applying the epsilon removal lemma in \cite[Corollary 7.5]{CLMP25}, it suffices to show: For any $0<\epsilon \ll 1$,
\begin{equation}\label{s2e11}
   D_p(R) \lesssim R^\epsilon
\end{equation}
holds for all $R\geq 1$ and $p>4$. Via a direct computation, we note
$$  |J_{\mathbf{Q}}(\xi;1)|\sim |\xi_2||\xi_3|.  $$
So by (\ref{s2s3e10}), we obtain
\begin{align}
D_p(R) \lesssim  &R^{\frac{\epsilon}{2}}\sup_{1<\mu_2,\mu_3\leq R^{1/2}}BD_p(1,\mu_2,\mu_3;R) +R^{\frac{\epsilon}{2}}D_p(1,R^{\frac{1}{2}},1;R)    +R^{\frac{\epsilon}{2}}D_p(1,1,R^{\frac{1}{2}};R) .  \label{s5c5e1a}
\end{align}
For the bilinear term, Theorem \ref{s3thm1} offers a nice bound. So, it suffices to consider the last two linear terms in (\ref{s5c5e1a}). Take the first linear term as an example. Since $\textbf{Q}$ satisfies the (CM) condition, we have the following sharp Stein--Tomas-type inequality (see \cite{Chr82,Moc96}):
$$   \|E^{\textbf{Q}}f\|_{L^{\frac{14}{3}}(B_R)} \lesssim \|f\|_{L^2}.  $$
On the other hand, we have the trivial $L^2$ estimate:
$$  \|E^{\textbf{Q}}f\|_{L^{2}(B_R)} \lesssim R\|f\|_{L^2}.    $$
By interpolation, we get 
$$   \|E^{\textbf{Q}}f\|_{L^{p}(B_R)} \lesssim R^{ \frac{7}{2p}-\frac{3}{4}  }\|f\|_{L^2}, $$
for $2\leq p \leq 14/3$. 
Thus, for any function $f$ with $\on{supp} f\subset [0,1]\times [b,b+R^{-1/2}]\times [0,1]\subset [0,1]^3$, we can apply H\"older's inequality to obtain
$$  \|E^{\textbf{Q}}f\|_{L^{p}(B_R)} \lesssim R^{ \frac{7}{2p}-\frac{3}{4}  } R^{-\frac{1}{2}(\frac{1}{2}-\frac{1}{p})}\|f\|_{L^p}=R^{\frac{4}{p}-1}\|f\|_{L^p}.   $$
Consequently, one gets
$$   R^{\frac{\epsilon}{2}}D_p(1,R^{\frac{1}{2}},1;R) \lesssim R^{\frac{4}{p}-1+\frac{\epsilon}{2}}\ll R^\epsilon,   $$
for $4\leq p \leq 14/3$. This completes the proof of (\ref{s2e11}).
\qed

\begin{remark}
    As a historical remark, this result has been obtained through the multilinear method and the $k$-linear method in {\rm\cite{GGO23,GO22,guo2023decoupling}}. 
\end{remark}

Finally, we consider more general quadratic forms $\mathbf{Q}$ with dimension $d$ and codimension $n$. Recall that in Case (2) of Theorem \ref{main thm}, $\mathbf{Q}$ consists of $(n-1)$ monomials and one polynomial, and we have some constraints on $P(\xi)$ and $\lambda_{2n-1},\lambda_{2n}$ to get sharp results. Here we provide more examples that are not confined to these limitations. It is worth pointing out that the proof requires more delicate techniques for the narrow part, compared with that of Theorem \ref{main thm}.

\begin{theorem}\label{s5t1}

\noindent $\rm{(1)}$     
    Suppose that 
    $$\mathbf{Q}(\xi)=(\xi_1^2,\xi_1\xi_2,...,\xi_1\xi_{s},\xi_\lambda \xi_{s+1}+P_{s+1}(\xi),...,\xi_{\lambda}\xi_n+P_n(\xi))$$
     is an $n$-tuple of quadratic forms on $\mathbb{R}^{n}$. Assume that $2\leq \lambda \leq s+1$.
Denote the total number of variables in $P_{s+1}(\xi),...,P_n(\xi)$ other than $\xi_1$ and $\xi_\lambda$ by $\theta$. We assume that one of the following two cases holds.
    \begin{itemize}
        \item[${\rm(a)}$] Each $P_j$ is independent of the variables $\xi_1,\xi_j,...,\xi_n$, and $s \geq \frac{n}{2}+\frac{\theta}{2}$.
        \item[${\rm(b)}$] There exists $P_j$ with no mixed terms that is independent of the variables $\xi_j,...,\xi_n$, $P_{j'}\equiv0$ for all $j'\neq j$, and $s \geq \frac{n}{2}+\frac{\theta}{4}$.
    \end{itemize}
    Then the estimate ${\rm(\ref{main aim})}$ holds for
    \begin{equation}\label{s5t1e1}
        q > s +3,\quad \frac{1}{p}+\frac{ s+2}{q}<1.
    \end{equation}

    \noindent $\rm{(2)}$      
    Let $k\geq 1$.	Suppose that 
    $$\mathbf{Q}(\xi)=(\xi_1^2,\xi_1\xi_2,...,\xi_1\xi_{n-1}, \xi_{1}\xi_n \pm \xi_{n+1}^2\pm\cdots\pm\xi_{n+k}^2  ).   $$
    Then the estimate ${\rm(\ref{main aim})}$ holds for
    \begin{equation}\label{s5t1e2}
        q > n +2,\quad \frac{1}{p}+\frac{ n+1}{q}<1.
    \end{equation}
    Moreover, the ranges ${\rm(\ref{s5t1e1})}$ and ${\rm(\ref{s5t1e2})}$ are all sharp up to the endpoints.
\end{theorem}

\noindent \textit{Proof sketch of Theorem \ref{s5t1}.}

(1) Note that in this case
$$   |J_{\mathbf{Q}}(\xi)| \sim |\xi_1|^{s}|\xi_\lambda|^{n-s}.   $$
By (\ref{s2s3e10}), we obtain
\begin{align}
D_p(R) \lesssim  &R^{\frac{\epsilon}{2}}\sup_{1<\mu_1,\mu_\lambda\leq R^{1/2}}BD_p(\mu_1,1,...,1,  \underset{\substack{\uparrow \\ \lambda}}{\mu_\lambda},1,...,1;R)  \nonumber \\
&+R^{\frac{\epsilon}{2}}D_p(R^{\frac{1}{2}},1,...,1;R)    +R^{\frac{\epsilon}{2}}D_p(1,...,1,\underset{\substack{\uparrow \\ \lambda}}{ R^{\frac{1}{2}}},1,...,1;R) .  \label{s5c5e1}
\end{align}
For the bilinear term, Theorem \ref{s3thm1} offers a nice bound by noting $s\geq n/2$. 

For the first linear term in (\ref{s5c5e1}), we can assume that each $P_j$ is independent of the variable $\xi_1$. This is because in Case (b), if $P_j$ depends on $\xi_1$, then its form must be $c\xi_1^2$, which can be dropped by adding a multiple of $\xi_1^2$ (the first component) to the $j$th component of $\mathbf{Q}$. Also note that $\lambda\neq 1$. Therefore, we can apply rescaling $\xi_1\rightarrow \xi_1/R^{1/2}$ and induction on scale to get
$$      R^{\frac{\epsilon}{2}} D_p(R^{\frac{1}{2}},1,...,1;R) \lesssim R^{\frac{s+3}{2p}-\frac{1}{2}+\frac{3}{2}\epsilon} \ll R^\epsilon$$
for $p>s+3$.

Next we consider the last linear term in (\ref{s5c5e1}). Suppose that we are in Case (a). Recall that each $P_j$ depends on $\theta$ variables other than $\xi_1$ and $\xi_\lambda$. We apply flat decoupling (\ref{d2}) to these variables to go from 1-scale to $R^{-1/2}$-scale. After this, by the locally constant property, the original $\mathbf{Q}$ can be reduced to  
$$ (\xi_1^2,\xi_1\xi_2,...,\xi_1\xi_{s},\xi_\lambda \xi_{s+1},...,\xi_{\lambda}\xi_n),  $$
which has been investigated in Case (1) in Theorem \ref{main thm}. Using rescaling and the relevant result, one concludes
\begin{equation*}
  R^{\frac{\epsilon}{2}}D_p(1,...,1,\underset{\substack{\uparrow \\ \lambda}}{ R^{\frac{1}{2}}},1,...,1;R)\lesssim R^{\frac{\theta}{2} (1-\frac{2}{p})+\frac{n+3\theta-s+3}{2p}-\frac{1+\theta}{2}+\epsilon}=R^{\frac{n+\theta-s+3}{2p}-\frac{1}{2}+\epsilon}\ll R^\epsilon
\end{equation*}
for $p> s+3$ due to $s \geq \frac{n}{2}+\frac{\theta}{2}$. So, the induction closes. 

Suppose that we are in Case (b).  Without loss of generality, we assume that $j=n$ and
    $$\mathbf{Q}(\xi)=(\xi_1^2,...,\xi_1\xi_{s},\xi_\lambda \xi_{s+1},...,\xi_\lambda \xi_{n-1},\xi_{\lambda}\xi_n\pm\xi_\alpha^2 \pm\cdots\pm\xi_\gamma^2),$$
with $2\leq \alpha,...,\gamma\leq n-1$, and $\alpha,...,\gamma\not=\lambda$. We apply flat decoupling (\ref{d2}) to the variables $\xi_\alpha,...,\xi_\gamma$ to go from $1$-scale to $R^{-1/4}$-scale, and then use the change of variables
\begin{equation*}
    \xi_\alpha\rightarrow\frac{\xi_\alpha}{R^{1/4}},\quad...\,,\quad \xi_\gamma\rightarrow\frac{\xi_\gamma}{R^{1/4}} ,   \quad  \xi_\lambda \rightarrow \frac{\xi_\lambda}{R^{1/2}}.
\end{equation*}
Notice that such rescaling keeps the surface $S_{\mathbf{Q}}$ unchanged. Covering $B_R$ with balls of scale $R/2$, and using induction on $R$, one concludes
\begin{equation*}
    R^{\frac{\epsilon}{2}}D_p(1,...,1,\underset{\substack{\uparrow \\ \lambda}}{ R^{\frac{1}{2}}},1,...,1;R) \lesssim R^{\frac{n+\theta/2-s+3}{2p}-\frac{1}{2}+\frac{3}{2}\epsilon} \ll R^\epsilon
\end{equation*}
for $p>s+3$ due to $s \geq \frac{n}{2}+\frac{\theta}{4}$. So, the induction closes.

(2) Take $i_1=n,...,i_k=n+k-1$, then
$$   |J_{\mathbf{Q}}(\xi;n,n+1,...,n+k-1)| \sim |\xi_1|^{n-1}|\xi_{n+k}|.   $$
By (\ref{s2s3e10}), we obtain
\begin{align}
D_p(R) \lesssim  &R^{\frac{\epsilon}{2}}\sup_{1<\mu_1,\mu_{n+k}\leq R^{1/2}}BD_p(\mu_1,1,...,1,\mu_{n+k};R)  \nonumber \\
&+R^{\frac{\epsilon}{2}}D_p(R^{\frac{1}{2}},1,...,1;R)    +R^{\frac{\epsilon}{2}} D_p(1,...,1,R^{\frac{1}{2}};R) .  \label{s5c5e2}
\end{align}
For the bilinear term, Theorem \ref{s3thm1} still offers a nice bound. 

We first consider the last linear term in (\ref{s5c5e2}). By the locally constant property, original $\mathbf{Q}$ can be reduced to 
    $$\mathbf{Q}_1(\xi)=(\xi_1^2,\xi_1\xi_2,...,\xi_{1}\xi_n \pm \xi_{n+1}^2\pm\cdots\pm\xi_{n+k-1}^2  ).   $$
We plan to use induction on the dimension $d$. When $d=n+1$ (base case), $\mathbf{Q}_1(\xi)=(\xi_1^2,...,\xi_{1}\xi_n)$. Let $g=g(\xi_1,...,\xi_n,x_{n+1})$ denote the inverse Fourier transform of $f$ in the variable $\xi_{n+1}$. By the result of Case (1) in Theorem \ref{main thm}, we know that 
$$    \|E^{\mathbf{Q}_1}g\|_{L^{n+3}(B_R^{2n})} \lesssim R^{\frac{\epsilon}{2}}\|g\|_{L^{n+3}(\mathbb{R}^{n})}.   $$
On the other hand, 
$$  \|E^{\mathbf{Q}_1}g\|_{L^{2}(B_R^{2n})} \lesssim R^{\frac{n}{2}}\|g\|_{L^2(\mathbb{R}^n)}.    $$
By interpolation, we have
$$   \|E^{\mathbf{Q}_1}g\|_{L^p(B_R^{2n})} \lesssim R^{n\cdot \frac{n+3}{n+1}(\frac{1}{p}-\frac{1}{n+3})+\frac{\epsilon}{2}}\|g\|_{L^p(\mathbb{R}^n)}, \quad \quad 2\leq p\leq  n+3.  $$
By this estimate for $\mathbf{Q}_1$ and the Hausdorff-Young inequality in the variable $\xi_{n+1}$, we get
\begin{align*}
  \|E^{\mathbf{Q}}f\|_{L^p(B_R^{2n+1})}  &\sim \|E^{\mathbf{Q}_1}f\|_{L^p(B_R^{2n+1})}  \\
  &  \leq  \Big\| \|E^{\mathbf{Q}_1}g\|_{L^p(B_R^{2n})} \Big\|_{L^p_{x_{n+1}}} \\   
  &  \lesssim  R^{n\cdot \frac{n+3}{n+1}(\frac{1}{p}-\frac{1}{n+3})+\frac{\epsilon}{2}}  \Big\|  \|g\|_{L^p_{\xi_1,...,\xi_n}}  \Big\|_{L^p_{x_{n+1}}}  \\
  &   \lesssim  R^{n\cdot \frac{n+3}{n+1}(\frac{1}{p}-\frac{1}{n+3})+\frac{\epsilon}{2}}  \Big\|  \|f\|_{L^{p'}_{\xi_{n+1}}}  \Big\|_{L^p_{ \xi_1,...,\xi_n}    }    \\
  &\lesssim R^{n\cdot \frac{n+3}{n+1}(\frac{1}{p}-\frac{1}{n+3})-\frac{1}{2}(\frac{1}{p'}-\frac{1}{p})+\frac{\epsilon}{2}} \|f\|_{L^p}    \\
  &=R^{\frac{n^2+4n+1}{(n+1)p} - \frac{3n+1}{2(n+1)}+\frac{\epsilon}{2}}\|f\|_{L^p}, 
\end{align*}
which implies that 
\begin{equation*}
 R^{\frac{\epsilon}{2}}D_p(1,...,1,R^{\frac{1}{2}};R) \lesssim R^{\frac{n^2+4n+1}{(n+1)p}-\frac{3n+1}{2(n+1)}+\frac{\epsilon}{2}} \ll R^{\epsilon},  \quad\quad  d=n+1,
\end{equation*}
for $p>n+2$. Having established the base case, for general cases $d=n+k$, we can use induction on the dimension $d$ to obtain
\begin{equation*}
 R^{\frac{\epsilon}{2}}D_p(1,...,1,R^{\frac{1}{2}};R) \ll R^{\epsilon},  \quad\quad  d\geq n+2,
\end{equation*}
for $p>n+2$.

Finally, we need to deal with the first linear term in (\ref{s5c5e2}). We will divide it into three subcases. First, we assume that $n=2$ and $k=1$. By using flat decoupling (\ref{d2}) in the variable $\xi_1$ to go from $R^{-1/2}$-scale to $R^{-1}$-scale, we get
$$  D_p(R^{\frac{1}{2}},1,1;R) \lesssim R^{\frac{1}{2}(1-\frac{2}{p})}D_p(R,1,1;R). $$
By the locally constant property, the original $\mathbf{Q}$ becomes $ \mathbf{Q}_2(\xi)=(0,\xi_3^2)$. Let $g=g(x_1,x_2,\xi_3)$ be the inverse Fourier transform of $f$ in the variables $\xi_1$ and $\xi_2$, and $\mathbf{Q}_2''(\xi)=(\xi_3^2)$. By the restriction estimate for the parabola \cite{Zyg74} and the Hausdorff--Young inequality in the variables $\xi_1$ and $\xi_2$, we have
\begin{align*}
    \|E^{\textbf{Q}}f\|_{L^p(B_R^{5})}& \sim  \|E^{\textbf{Q}_2}f\|_{L^p(B_R^{5})} \\
    &\lesssim  R^{\frac{1}{p}} \Big\| \|E^{\textbf{Q}_2''}g\|_{L^p(B_R^{2})}  \Big\|_{L^p_{x_1,x_2}}   \\
    & \lesssim R^{\frac{1}{p}} \Big\| \|g\|_{L^p_{\xi_3}} \Big\|_{L^p_{x_1,x_2}}  \\
    & \lesssim R^{\frac{1}{p}}  \Big\| \|f\|_{L^{p'}_{\xi_1,\xi_2}} \Big\|_{L^p_{\xi_{3}}} \\
    & \lesssim R^{\frac{1}{p}-(\frac{1}{p'}-\frac{1}{p})}\|f\|_{L^p},
\end{align*}
and so
$$    R^{\frac{\epsilon}{2}}D_p(R^{\frac{1}{2}},1,1;R) \lesssim R^{\frac{1}{2}(1-\frac{2}{p})+\frac{1}{p}-(\frac{1}{p'}-\frac{1}{p})+\frac{\epsilon}{2}} =R^{\frac{2}{p}-\frac{1}{2}+\frac{\epsilon}{2}}\ll R^\epsilon $$
for $p > n+2=4$. Next, we assume that $n = 3,4$  and $k=1$. We use flat decoupling (\ref{d2}) in the variable $\xi_n$ to go from $1$-scale to $R^{-1/2}$-scale, then
$$ D_p(R^{\frac{1}{2}},1,...,1;R)\lesssim R^{\frac{1}{2}(1-\frac{2}{p})} D_p(R^{\frac{1}{2}},1,...,1,R^{1/2},1;R) . $$
In this setting, by the locally constant property, the original $\mathbf{Q}$ is reduced to 
$$   \mathbf{Q}_3(\xi)=(\xi_1^2,\xi_1\xi_2,...,\xi_1\xi_{n-1},\xi_{n+1}^2).  $$
Using the result of Case (1) in Theorem \ref{main thm}, we have
  $$    \|E^{\mathbf{Q}_3}g\|_{L^{n+2}(B_R^{2n})} \lesssim R^{\frac{\epsilon}{2}}\|g\|_{L^{n+2}(\mathbb{R}^{n})}.   $$   
By this estimate for $\mathbf{Q}_3$ and the Hausdorff--Young inequality in the variable $\xi_{n}$ as before, we get
$$    R^{\frac{\epsilon}{2}}  D_p(R^{\frac{1}{2}},1,...,1;R)\lesssim  R^{\frac{1}{2}(1-\frac{2}{p})    -\frac{1}{2}    (\frac{1}{p'}-\frac{1}{p})+\epsilon}\ll R^\epsilon, $$
for $p>n+2$.
For all other $n$ and $k$, we use decoupling for the hyperbolic paraboloids (\ref{d10}) in the variables $\xi_{n+1},...,\xi_{n+k}$ to go from 1-scale to $R^{-1/4}$-scale, followed by rescaling and induction on $R$, to obtain
\begin{align*}
    R^{\frac{\epsilon}{2}}D_p(R^{\frac{1}{2}},1,...,1;R) &\lesssim R^{\frac{k}{4}-\frac{k+1}{2p}+\frac{\epsilon}{2}}  D_p(R^{\frac{1}{2}},1,...,1,R^{\frac{1}{4}},...,R^{\frac{1}{4}};R)   \\
    &\leq R^{\frac{k}{4}-\frac{k+1}{2p}+\frac{n+k+3}{2p}-\frac{k+2}{4}+\frac{\epsilon}{2}}  D_p(R) \\
    & \lesssim   R^{\frac{n+2}{2p}-\frac{1}{2}+\frac{3}{2}\epsilon}\\
    &\ll R^\epsilon
\end{align*}
for $p>n+2$. So, the induction closes. 

\qed

From the above arguments, it is plausible that our method is powerful in the study for degenerate higher codimensional surfaces. Nevertheless, it does not cover all possible cases, even when the quadratic surface is monomial. 
For example, as we have pointed out, when
$$\mathbf{Q}(\xi)=(\xi_{1}\xi_2,\xi_{2}\xi_3,\xi_{3}\xi_4,\xi_{4}\xi_1),$$
we can compute $J_{\mathbf{Q}}(\xi)\equiv 0$, and then the condition in Theorem \ref{s3thm1} is not met at all. Moreover, even if our method can be applied, it may well be the case that the restriction estimate we obtain is not sharp. A typical example is 
$$\mathbf{Q}(\xi)=(\xi_1^2,...,\xi_1 \xi_{n-1},\xi_{1}\xi_n+\xi_{2}^2).$$
Through the same argument as in the last part of the previous section, we expect that (\ref{main aim}) may hold for $p>n+2$ if $p=q$. Our method only gives the range $p>n+3$, whose main loss comes from Theorem \ref{s3thm1}, which only covers $p>n+3$. Therefore, our method cannot yield a sharp bound of here. If we consider the classical broad--narrow analysis, then the constraint of the broad part is $p>4$, while the narrow part requires $p>n+7/2$ due to the loss of decoupling.

\vskip0.3cm

\section*{Acknowledgements} 
The authors would like to thank the anonymous referee for valuable comments and suggestions that improved the clarity and presentation of the manuscript. This work was supported by the National Key R\&D program of China (Grant number 2022YFA 1005700). Zhenbin Cao was supported by High-level Talent Research Start-up Project Funding of Henan Academy of Sciences (Grant number 241819032); and the National Natural Science Foundation of China (Grant numbers 12401117, 12371100). Changxing Miao was supported by the National Natural Science Foundation of China (Grant numbers 12531005, 12371095).

\section*{Data availability statement}
Data sharing is not applicable to this article.


\end{document}